\input amstex\documentstyle {amsppt}  
\pagewidth{12.5 cm}\pageheight{19 cm}\magnification\magstep1
\topmatter
\title Unipotent elements in small characteristic, III\endtitle
\author G. Lusztig\endauthor
\address Department of Mathematics, M.I.T., Cambridge, MA 02139\endaddress
\thanks Supported in part by the National Science Foundation\endthanks
\endtopmatter   
\document
\define\gr{\text{\rm gr}}

\define\dy{\dot y}

\define\be{\bar e}

\define\bT{\bar T}
\define\bQ{\bar Q}

\define\frl{\forall}

\define\si{\sim}

\define\sqc{\sqcup}

\define\qua{\quad}

\define\dx{\dot x}

\define\bA{\bar A}
\define\baf{\bar f}

\define\lb{\linebreak}

\define\op{\oplus}

\define\part{\partial}
\define\em{\emptyset}

\define\ra{\rangle}
\define\n{\notin}

\define\m{\mapsto}
\define\do{\dots}
\define\la{\langle}
\define\bsl{\backslash}

\define\lra{\leftrightarrow}

\define\sm{\smallmatrix}
\define\esm{\endsmallmatrix}
\define\sub{\subset}    

\define\T{\times}
\define\ti{\tilde}
\define\nl{\newline}
\redefine\i{^{-1}}

\define\Ad{\text{\rm Ad}}
\define\Hom{\text{\rm Hom}}
\define\End{\text{\rm End}}

\define\a{\alpha}
\redefine\b{\beta}

\define\g{\gamma}
\redefine\d{\delta}

\define\et{\eta}

\redefine\o{\omega}

\define\ph{\phi}
\define\ps{\psi}

\define\x{\xi}

\redefine\D{\Delta}

\define\Si{\Sigma}
\define\Th{\Theta}

\define\Ph{\Phi}
\define\Ps{\Psi}

\define\kk{\bold k}

\define\CC{\bold C}

\define\FF{\bold F}

\define\NN{\bold N}

\define\TT{\bold T}

\define\WW{\bold W}
\define\ZZ{\bold Z}

\define\cm{\Cal M}
\define\cn{\Cal N}
\define\co{\Cal O}

\define\car{\Cal R}

\define\cu{\Cal U}

\define\fg{\frak g}

\define\fo{\frak o}

\define\fs{\frak s}

\define\fA{\frak A}

\define\fD{\frak D}

\define\fF{\frak F}

\define\fS{\frak S}

\define\tm{\ti m}

\define\tG{\ti G}

\define\tV{\ti V}

\define\KO{K}
\define\LN{L1}
\define\LUU{L2}
\define\LU{L3}
\define\SPR{S}
\define\ST{St}

\head 0. Introduction\endhead
Let $\kk$ be an algebraically closed field of characteristic exponent $p\ge1$. Let $G$ be a connected reductive 
algebraic group over $\kk$ and let $\fg$ be the Lie algebra of $G$. Note that $G$ acts on $G$ and on $\fg$ by the
adjoint action. Let $\cu_G$ be the variety of unipotent elements of $G$. Let $\cn_\fg$ be the variety of nilpotent
elements of $\fg$. In \cite{\LUU}, we have proposed a definition of a partition of $\cu_G$ into smooth locally 
closed $G$-stable pieces which are indexed by the unipotent classes in the group over $\CC$ of the same type as 
$G$ and which in many ways seem to depend very smoothly on $p$; moreover we studied in detail the pieces of 
$\cu_G$ for types $A$ and $C$. In \cite{\LU} we have studied in detail the pieces of $\cu_G$ for types $B$ and 
$D$; however the definition in \cite{\LU} was not based on the proposal of \cite{\LUU} (which involved the 
partial order of unipotent classes). 

In this paper we propose another general definition of the pieces of $\cu_G$ which is not 
based on the partial order of unipotent classes and we show that this new definition unifies the definitions in 
\cite{\LUU}, \cite{\LU} in the sense that for types $A,C$ it can be identified with the definition in \cite{\LUU}
while for types $B,D$ it can be identified with the definition in \cite{\LU}. The idea of the new definition is as
follows. One needs to consider a grading $\fg=\op_n\fg_n$ analogous to the one associated to a nilpotent element 
over $\CC$ by the Morozov-Jacobson theorem. Let $G_0$ be the closed connected subgroup of $G$ corresponding to 
$\fg_0$. Now $G_0$ acts naturally on $\fg_2$. The main ingredient in the definition of a piece is the definition 
of a suitable open $G_0$-invariant subset $\fg_2^!$ of $\fg_2$. When $p=1$ or $p\gg0$, the subset $\fg_2^!$ is by
definition the unique open $G_0$-orbit in $\fg_2$. But this definition is not correct for general $p$. Now note 
that when $p=1$ the centralizer in $G$ of any element of the open $G_0$-orbit in $\fg_2$ is contained in 
$G_{\ge0}$, the parabolic subgroup of $G$ whose Lie algebra is $\sum_{i\ge0}\fg_i$ (a result of Kostant 
\cite{\KO}). Based on this, we propose to define (for general $p$) the set $\fg_2^!$ as the set of all $x\in\fg_2$
such that the centralizer of $x$ in $G$ is contained in $G_{\ge0}$. It turns out that (at least in types 
$A,B,C,D$) this condition gives exactly the unique open $G_0$-orbit in $\fg_2$ when $p=1$ or $p\gg0$, while in 
general it defines a subset of $\fg_2$ which is a union of possibly several $G_0$-orbits (their number is a power
of $2$) but which is exactly what is needed to define the pieces of $\cu_G$.

This paper is organized as follows. In Section 1 we give the definition of the sets $\fg_2^!$ and show that in 
types $A,B,C,D$ these sets can be identified with certain explicit subsets of $\fg_2$ considered in \cite{\LUU},
\cite{\LU}. In Section 2 we use the sets $\fg_2^!$ to define the pieces of $\cu_G$. In the appendix (written by 
the author and T. Xue) we define (at least in types $A,B,C,D$) a partition of $\cn_\fg$ into smooth locally closed
$G$-stable pieces which are indexed by the unipotent classes in the group over $\CC$ of the same type as $G$ 
(using again the subsets $\fg_2^!$ defined in Section 1).

{\it Notation.} The cardinal of a finite set $X$ is denoted by $|X|$.

{\it Errata to \cite{\LN}.} 

On p.206, the last sentence in 6.8, "For classical types...classes", should be removed.

{\it Errata to \cite{\LUU}.} 

On p.452 line -5 replace $\Pi\cup\Th$ by $\Th$.

On p.452 line -4 replace $\Pi\cup\ti\Th$ by $\ti\Th$.

{\it Errata to \cite{\LU}.} 

On p.774 line 6: replace "an (injective)" by "a".

On p.778, line 13: replace last $T^a$ by $T^ax$.

On p.779, line 10: after "(b) and (c)" insert: "and $f_a=f_{-a}$ for all $a$".

On p.797, remove last line of 4.2 and last line of 4.3.

\head Contents\endhead
1. The sets $\fg_2^{\d!}$.

2. The pieces in the unipotent variety of $G$.

Appendix (by G. Lusztig and T. Xue): The pieces in the nilpotent variety of $\fg$.

\head 1. The sets $\fg_2^{\d!}$\endhead
\subhead 1.1\endsubhead
Let $\TT,\WW$ be "the maximal torus" and "the Weyl group" of $G$ and let $Y_G=\Hom(\kk^*,\TT)$. Note that $\WW$ 
acts naturally on $\TT$ and on $Y_G$.

Let $G'$ be a connected reductive algebraic group over $\CC$ of the same type as $G$. In particular we have 
canonically $Y_G=Y_{G'}$ compatibly with the $\WW$-actions. Note that $G$ acts by conjugation on $\Hom(\kk^*,G)$;
similarly, $G'$ acts by conjugation on $\Hom(\CC^*,G')$. We have canonically $G\bsl\Hom(\kk^*,G)=\WW\bsl Y_G$, 
$G'\bsl\Hom(\CC^*,G')=\WW\bsl Y_{G'}$. Let $\fD_{G'}$ be the set of all $f\in\Hom(\CC^*,G')$ such that for some 
homomorphism of algebraic groups $h:SL_2(\CC)@>>>G'$ we have $h\left(\sm a&0\\0&a\i\esm\right)=f(a)$ for all 
$a\in\CC^*$. Let $\fD_G$ be the set of all $\d\in\Hom(\kk^*,G)$ such that the image of $\d$ in 
$G\bsl\Hom(\kk^*,G)=\WW\bsl Y_G=\WW\bsl Y_{G'}=G'\bsl\Hom(\CC^*,G')$ can be represented by an element in 
$\fD_{G'}\sub\Hom(\CC^*,G')$. (See \cite{\LUU, 1.1}.) Note that $\fD_G$ is a union of $G$-orbits. The $G'$-orbits
in $\fD_{G'}$ were classified by Dynkin. They form a finite set in natural bijection  with the set of $G$-orbits 
in $\fD_G$.

Let $G^{der}$ be the derived group of $G$. The simply connected covering $\tG^{der}@>>>G$ of $G^{der}$ induces a 
bijection $\fD_{\tG^{der}}@>>>\fD_G$. This follows from the analogous assertion for $G'$ which is immediate.

\subhead 1.2\endsubhead
Let $\d\in\Hom(\kk^*,G)$. For any $i\in\ZZ$ we set

$\fg_i^\d=\{x\in\fg;\Ad(\d(a))x=a^ix\qua\frl a\in\kk^*\}$.
\nl
We have a direct sum decomposition $\fg=\op_{i\in\ZZ}\fg_i^\d$. For $i\in\NN$ we set 
$\fg_{\ge i}^\d=\op_{j\in\ZZ;j\ge i}\fg_j^\d$; note that $\fg_{\ge i}^\d$ is the Lie algebra of a well defined 
closed connected subgroup $G_{\ge i}^\d$ of $G$. We have $\do\sub G_{\ge2}^\d\sub G_{\ge1}^\d\sub G_{\ge0}^\d$ and
$G_{\ge0}^\d$ is a parabolic subgroup of $G$ with unipotent radical $G_{\ge1}^\d$ and with Levi subgroup $G_0^\d$
(with Lie algebra $\fg_0^\d$); moreover, for any $i$, $G_{\ge i}^\d$ is a normal subgroup of $G_{\ge0}^\d$.

We set $\fg_{<0}^\d=\op_{j\in\ZZ;j<0}\fg_j^\d$; note that $\fg_{<0}^\d$ is the Lie algebra of a well defined 
closed connected subgroup $G_{<0}^\d$ of $G$. We have $G_{<0}^\d\cap G_{\ge0}^\d=\{1\}$.

For any $x\in\fg$ let $G_x=\{g\in G;\Ad(g)x=x\}$. Let
$$\fg_2^{\d!}=\{x\in\fg_2^\d;G_x\sub G^\d_{\ge0}\}.$$
For $h\in G^\d_0$, $x\in\fg$ we have $G_{\Ad(h)x}=hG_xh\i$; hence $\Ad(h)\fg_2^{\d!}=\fg_2^{\d!}$. Thus, 
$\fg_2^{\d!}$ is a union of orbits for the $\Ad$-action of $G^\d_0$ on $\fg^\d_2$.

In 1.3, 1.4, 1.5 we will describe explicitly the set $\fg_2^{\d!}$ in a number of cases.

\subhead 1.3\endsubhead
In this subsection we fix a $\kk$-vector space $V$ of finite dimension and we assume that $G=GL(V)$. Then
$\fg=\End(V)$. A $\ZZ$-grading $V=\op_{i\in\ZZ}V_i$ of $V$ is said to be {\it good} if 
$\dim V_i=\dim V_{-i}\ge\dim V_{-i-2}$ for all $i\ge0$. If a good $\ZZ$-grading $(V_i)$ is given and $a\in\ZZ$, we
set

$V_{\ge a}=\sum_{j;j\ge a}V_a$,

$\End(V)_a=\{A\in\fg;A(V_r)\sub V_{r+a}\qua\frl r\in\ZZ\}$,

$\End(V)_{\ge a}=\{A\in\fg;A(V_{\ge r})\sub V_{\ge r+a}\qua\frl r\in\ZZ\}$.

To give an element $\d\in\fD_G$ is the same as to give a good grading $V=\op_{i\in\ZZ}V_i$ of $V$ ($\d$ is given 
in terms of the $\ZZ$-grading by $\d(a)|_{V_r}=a^r$ for all $a\in\kk^*,r\in\ZZ$.) 

In the remainder of this subsection we fix $\d\in\fD_G$ and let $(V_i)$ be the corresponding good grading of $V$.
Then $\fg_i^\d=\End(V)_i$. Let 
$\End(V)_2^0$ be the set of all $A\in\End(V)_2$ such that $A^n:V_{-n}@>>>V_n$ is an isomorphism for any $n\ge0$.
It is easy to see that $\End(V)_2^0\ne\em$. The following result will be proved in 1.6.

(a) $\fg_2^{\d!}=\End(V)_2^0$. 

\subhead 1.4\endsubhead
In this subsection we fix a $\kk$-vector space $V$ of finite even dimension with a fixed nondegenerate symplectic
form $(,):V\T V@>>>\kk$ and we assume that $G=Sp(V)=\{T\in GL(V);T\text{ preserves }(,)\}$. Let 
$\fs(V)=\{T\in\End(V);(Tv,v')+(v,Tv')=0\qua\frl v,v'\in V\}$. Then $\fg=\fs(V)$.

A $\ZZ$-grading $V=\op_{i\in\ZZ}V_i$ of $V$ is said to be {\it $s$-good} if it is good (see 1.3), $\dim V_i$ is 
even for any even $i$ and $(V_i,V_j)=0$ whenever $i+j\ne0$.

To give an element $\d\in\fD_G$ is the same as to give an $s$-good grading $V=\op_{i\in\ZZ}V_i$ of $V$ ($\d$ is 
given in terms of the $\ZZ$-grading by $\d(a)|_{V_r}=a^r$ for all $a\in\kk^*,r\in\ZZ$.) 

In the remainder of this subsection we fix $\d\in\fD_G$ and let $(V_i)$ be the corresponding $s$-good grading of 
$V$. Let $\fs(V)_i=\fs(V)\cap\End(V)_i$. Then $\fg_i^\d=\fs(V)_i$ for any $i$. Let 
$\fs(V)_2^0=\fs(V)_2\cap\End(V)_2^0$ (notation of 1.3). The following result will be proved in 1.7.

(a) $\fg_2^{\d!}=\fs(V)_2^0$. 

\subhead 1.5\endsubhead
In this subsection we fix a $\kk$-vector space $V$ of finite dimension with a fixed nondegenerate quadratic form
$Q:V@>>>\kk$ with associate symmetric bilinear form $(,):V\T V@>>>\kk$. (Recall that 
$(v,v')=Q(v+v')-Q(v)-Q(v')$ for $v,v'\in V$ and that the nondegeneracy of $Q$ means that, if $\car$ is the 
radical of $(,)$, then $\car=0$ if $p\ne2$ and $Q:\car@>>>\kk$ is injective if $p=2$.) We assume that $G=SO(V)$, 
the identity component of $O(V)=\{T\in GL(V);Q(T(v))=Q(v)\qua\frl v\in V\}$. Let 
$\fo(V)=\{T\in\End(V);(Tv,v)=0\qua\frl v\in V, T|_{\car}=0\}$. Then $\fg=\fo(V)$.

A $\ZZ$-grading $V=\op_{i\in\ZZ}V_i$ of $V$ is said to be {\it $o$-good} if it is good (see 1.3), $\dim V_i$ is 
even for any odd $i$, $(V_i,V_j)=0$ whenever $i+j\ne0$ and $Q|_{V_i}=0$ whenever $i\ne0$.

To give an element $\d\in\fD_G$ is the same as to give an $o$-good grading $V=\op_{i\in\ZZ}V_i$ of $V$ ($\d$ is 
given in terms of the $\ZZ$-grading by $\d(a)|_{V_r}=a^r$ for all $a\in\kk^*,r\in\ZZ$.) 

In the remainder of this subsection we fix $\d\in\fD_G$ and let $(V_i)$ be the corresponding $o$-good grading of 
$V$. Let $\fo(V)_i=\fo(V)\cap\End(V)_i$. Then $\fg_i^\d=\fo(V)_i$ for any $i$. Let $\fo(V)_2^0$ be the set of all 
$A\in\fo(V)_2$ such that

(i) for any odd $n\ge1$, $A^n:V_{-n}@>>>V_n$ is an isomorphism;

(ii) for any even $n\ge0$, $A^{n/2}:V_{-n}@>>>V_0$ is injective and the restriction of $Q$ to $A^{n/2}(V_{-n})$ is
nondegenerate.
\nl
Note that the equality $\fo(V)_2^0=\fo(V)_2\cap\End(V)_2^0$ (notation of 1.3) holds when $p\ne2$, but not 
necessarily when $p=2$. The following result will be proved in 1.8.

(a) $\fg_2^{\d!}=\fo(V)_2^0$. 

\subhead 1.6\endsubhead
We prove 1.3(a). Generally, $x_k$ will denote an element of $V_k$. Let $A\in\End(V)_2-\End(V)_2^0$. The transpose
$A^*:V^*@>>>V^*$ of $A$ carries $V_i^*$ to $V_{i-2}^*$. 

Assume first that $A:V_{-i}@>>>V_{-i+2}$ is not injective for some $i\ge2$. Then $A^*:V_{-i+2}^*@>>>V_{-i}^*$ is 
not surjective. Since $\dim V_{-i+2}\ge\dim V_{-i}$ it follows that $A^*:V_{-i+2}^*@>>>V_{-i}^*$ is not injective.
We can find $e_{-i}\in V_{-i}-\{0\}$ such that $Ae_{-i}=0$ and $\x_{-i+2}\in V_{-i+2}^*-\{0\}$ such that 
$A^*\x_{-i+2}=0$. Define $B\in\End(V)$ by 

$B(\sum_kx_k)=\sum_{k\ne-i}x_k+(x_{-i}+\x_{-i+2}(x_{-i+2})e_{-i})$.
\nl
We have $B\in G^\d_{<0}-\{1\}$. Hence $B\n G_{\ge0}^\d$. We have    
$$\align&BA(\sum_kx_k)-AB(\sum_kx_k)=\sum_{k\ne-i}Ax_{k-2}+
(Ax_{-i-2}+\x_{-i+2}(Ax_{-i})e_{-i})\\&-\sum_{k\ne-i}x_k-(Ax_{-i}+\x_{-i+2}(x_{-i+2})Ae_{-i})\\&
=(A^*\x_{-i+2})(x_{-i})e_{-i}-\x_{-i+2}(x_{-i+2})Ae_{-i}=0.\endalign$$
Thus $AB=BA$. We see that $A\n\fg^{\d!}_2$.

Assume next that $A:V_i@>>>V_{i+2}$ is not surjective for some $i\ge0$. Since $\dim V_i\ge\dim V_{i+2}$, 
$A:V_i@>>>V_{i+2}$ is not injective. Moreover $A^*:V_{i+2}^*@>>>V_i^*$ is not injective. We can find 
$e_i\in V_i-\{0\}$ such that $Ae_i=0$ and $\x_{i+2}\in V_{i+2}^*-\{0\}$ such that $A^*\x_{i+2}=0$. Define 
$B\in\End(V)$ by 

$B(\sum_kx_k)=\sum_{k\ne i}x_k+(x_i+\x_{i+2}(x_{i+2})e_i)$. 
\nl
We have $B\in G^\d_{<0}-\{1\}$. Hence $B\n G_{\ge0}^\d$. We have 
$$\align&BA(\sum_kx_k)-AB(\sum_kx_k)=\sum_{k\ne i}Ax_{k-2}+(Ax_{i-2}+\x_{i+2}(Ax_i)e_i)\\&
-\sum_{k\ne i}x_k-(Ax_i+\x_{i+2}(x_{i+2})Ae_i)=(A^*\x_{i+2})(x_i)e_i-\x_{i+2}(x_{i+2})Ae_i=0.\endalign$$
Thus $AB=BA$. We see that $A\n\fg^{\d!}_2$.

We now assume that $A:V_{-i}@>>>V_{-i+2}$ is injective for any $i\ge2$ and $A:V_i@>>>V_{i+2}$ is surjective for 
any $i\ge0$. For some $n>0$, $A^n:V_{-n}@>>>V_{n}$ is not an isomorphism hence $A^{*n}:V_n^*@>>>V_{-n}^*$ is not 
an isomorphism. We can find $e_{-n}\in V_{-n}-\{0\}$ such that $A^ne_{-n}=0$. We can find $\x_n\in V_n^*-\{0\}$ 
such that $A^{*n}\x_n=0$. For any $j\ge0$ we set $e_{2j-n}=A^je_{-n}\in V_{2j-n}$,
$\x_{n-2j}=A^{*j}\x_n\in V_{n-2j}^*$. Note that $e_n=0,\x_{-n}=0$. Also, $e_m\ne0$ if $m\le0$, $m=n\mod2$ and 
$\x_m\ne0$ if $m\ge0$, $m=n\mod2$. Define $B\in\End(V)$ by
$$B(\sum_kx_k)=\sum_{k\n\{2h-n;h\in[0,n-1]\}}x_k
+\sum_{j\in[0,n-1]}(x_{2j-n}+\x_{2j-n+2}(x_{2j-n+2})e_{2j-n}).$$
We have $B\in G^\d_{<0}$. If $n$ is even then the term corresponding to $j=n/2-1$ is $x_{-2}+\x_0(x_0)e_{-2}$
and $\sum_kx_k\m\x_0(x_0)e_{-2}$ is $\ne0$ since $e_{-2}\ne0,\x_0\ne0$. If $n$ is odd then the term corresponding
to $j=(n-1)/2$ is $x_{-1}+\x_1(x_1)e_{-1}$ and $\sum_kx_k\m\x_1(x_1)e_{-1}$ is $\ne0$ since $e_{-1}\ne0,\x_1\ne0$.
Thus $B\ne1$ so that $B\n G_{\ge0}^\d$. We have
$$\align&BA(\sum_kx_k)-AB(\sum_kx_k)\\&=
\sum_{k\n\{2h-n;h\in[0,n-1]\}}Ax_{k-2}+\sum_{j\in[0,n-1]}(Ax_{2j-n-2}+\x_{2j-n+2}(Ax_{2j-n})e_{2j-n})\\&
-\sum_{k\n\{2h-n;h\in[0,n-1]\}}Ax_k-\sum_{j\in[0,n-1]}(Ax_{2j-n}+\x_{2j-n+2}(x_{2j-n+2})Ae_{2j-n})\\&
=\sum_{j\in[0,n-1]}(A^*\x_{2j-n+2})(x_{2j-n})e_{2j-n}-\sum_{j\in[0,n-1]}\x_{2j-n+2}(x_{2j-n+2})Ae_{2j-n}\\&
=\sum_{j\in[0,n-1]}\x_{2j-n}(x_{2j-n})e_{2j-n}-\sum_{j\in[0,n-1]}\x_{2j-n+2}(x_{2j-n+2})e_{2j-n+2}\\&
=\sum_{j\in[0,n-1]}\x_{2j-n}(x_{2j-n})e_{2j-n}-\sum_{j\in[1,n]}\x_{2j-n}(x_{2j-n})e_{2j-n}\\&
=\x_{-n}(x_{-n})e_{-n}-\x_n(x_n)e_n=0\endalign$$
since $x_{-n}=0,e_n=0$. Thus $BA=AB$. We see that $A\n\fg^{\d!}_2$.

We have shown that $\End(V)_2-\End(V)_2^0\sub\End(V)_2-\fg^{\d!}_2$.

Conversely, let $A\in\End(V)_2^0$. We show that $A\in\fg^{\d!}_2$. Let $B\in G$ be such that $AB=BA$. It is enough
to show that $B\in G_{\ge0}^\d$. We argue by induction on $\dim V$.
If $V=0$ the result is clear. Assume now that $V\ne0$. Let $m$ be the largest integer $\ge0$ such that
$V_m\ne0$. If $m=0$ we have $G_{\ge0}^\d=G$ and the result is clear. Assume now that $m\ge1$. We have $A^mV=V_m$,
$\ker(A^m:V@>>>V)=V_{\ge-m+1}$. Since $BA=AB$ we have $B(A^mV)=A^mV$, $B(\ker(A^m:V@>>>V))=\ker(A^m:V@>>>V)$. 
Hence $B(V_m)=V_m$ and $B(V_{\ge-m+1})=V_{\ge-m+1}$. Hence $B$ induces an automorphism $B':V'@>>>V'$ where 
$V'=V_{\ge-m+1}/V_m$. We have canonically $V'=V_{-m+1}\op V_{-m+2}\op\do\op V_{m-1}$ and $\End(V')_2,\End(V')_2^0$
are defined in terms of this (good) grading. Now $A$ induces an element $A'\in\End(V')_2^0$ and we
have $B'A'=A'B'$. By the induction hypothesis, for any $i\in[-m+1,m-1]$, the subspace $V_i+V_{i+1}+\do+V_{m-1}$ of
$V'$ is $B'$-stable. Hence the subspace $V_i+V_{i+1}+\do+V_{m-1}+V_m$ of $V$ is $B$-stable. We see that
$B\in G_{\ge0}^\d$. This completes the proof of 1.3(a).

\subhead 1.7\endsubhead
We prove 1.4(a). Generally, $x_k$ will denote an element of $V_k$. 

Let $A\in\fs(V)_2-\fs(V)_2^0$. Assume first that $A:V_{-i}@>>>V_{-i+2}$ is not injective for some $i>2$. Then 
$A:V_{i-2}@>>>V_i$ is not surjective and since $\dim V_{i-2}\ge\dim V_i$, we see that $A:V_{i-2}@>>>V_i$ is not 
injective. We can find $e_{-i}\in V_{-i}-\{0\}$ such that $Ae_{-i}=0$. We can find $e_{i-2}\in V_{i-2}-\{0\}$ such
that $Ae_{i-2}=0$. Define $B\in\End(V)$ by 
$$B(\sum_kx_k)=\sum_{k\ne-i,i-2}x_k+(x_{-i}+(e_{i-2},x_{-i+2})e_{-i})+(x_{i-2}+(e_{-i},x_i)e_{i-2})$$
where $x_k\in V_k$. We have
$$\align&(B(\sum_kx_k),B(\sum_kx'_k))-(\sum_kx_k,\sum_kx'_k)\\&=
(x_{-i}+(e_{i-2},x_{-i+2})e_{-i},x'_i)+(x_i,x'_{-i}+(e_{i-2},x'_{-i+2})e_{-i})\\&
+(x_{i-2}+(e_{-i},x_i)e_{i-2},x'_{-i+2})+(x_{-i+2},x'_{i-2}+(e_{-i},x'_i)e_{i-2})\\&
+\sum_{k\ne-i+2,-i,i,i-2}(x_{-k},x'_k)-\sum_k(x_{-k},x'_k)\\&
=(e_{i-2},x_{-i+2})(e_{-i},x'_i)+(x_i,e_{-i})(e_{i-2},x'_{-i+2})\\&
+(e_{-i},x_i)(e_{i-2},x'_{-i+2})+(x_{-i+2},e_{i-2})(e_{-i},x'_i)\\&=0.\endalign$$
Thus, $B\in Sp(V)$. More precisely, $B\in G^\d_{<0}-\{1\}$. Hence $B\n G_{\ge0}^\d$. We have 
$$\align&BA(\sum_kx_k)-AB(\sum_kx_k)\\&=
\sum_{k\ne-i,i-2}Ax_{k-2}+(Ax_{-i-2}+(e_{i-2},Ax_{-i})e_{-i})+(Ax_{i-4}+(e_{-i},Ax_{i-2})e_{i-2})\\&
-\sum_{k\ne-i,i-2}Ax_k-(Ax_{-i}+(e_{i-2},x_{-i+2})Ae_{-i})-(Ax_{i-2}+(e_{-i},x_i)Ae_{i-2})\\&
=-(Ae_{i-2},x_{-i})e_{-i}-(Ae_{-i},x_{i-2})e_{i-2}-(e_{i-2},x_{-i+2})Ae_{-i}-(e_{-i},x_i)Ae_{i-2}=0.\endalign$$
Thus, $AB=BA$. We see that $A\n\fg^{\d!}_2$.

Next we assume that $A:V_{-2}@>>>V_0$ is not injective. We can find $e_{-2}\in V_{-2}-\{0\}$ such that 
$Ae_{-2}=0$. Note that $K'=\ker(A:V_0@>>>V_2)\ne0$. (Indeed $A:V_0@>>>V_2$ is the transpose of $A:V_{-2}@>>>V_0$
hence is not surjective. But $\dim V_0\ge\dim V_2$ hence $K'\ne0$.) We can find $e_0\in V_0-\{0\}$ such that 
$Ae_0=0$. Define $B\in\End(V)$ by 
$$B(\sum_kx_k)=\sum_{k\ne-2,0}x_k+(x_{-2}+(e_0,x_0)e_{-2})+(x_0+(e_{-2},x_2)e_0).$$
We have
$$\align&(B(\sum_kx_k),B(\sum_kx'_k))-(\sum_kx_k,\sum_kx'_k)\\&
=(x_0+(e_{-2},x_2)e_0,x'_0+(e_{-2},x'_2)e_0)+(x_{-2}+(e_0,x_0)e_{-2},x'_2)\\&+(x_2,x'_{-2}+(e_0,x'_0)e_{-2})+
\sum_{k\ne-2,0,2}(x_{-k},x'_k)+\sum_k(x_{-k},x'_k)\\&
=(x_0,(e_{-2},x'_2)e_0)+((e_{-2},x_2)e_0,x'_0)+((e_0,x_0)e_{-2},x'_2)+(x_2,(e_0,x'_0)e_{-2})\\&
=(x_0,e_0)(e_{-2},x'_2)+(e_{-2},x_2)(e_0,x'_0)+(e_0,x_0)(e_{-2},x'_2)+(x_2,e_{-2})(e_0,x'_0)=0.\endalign$$
Thus $B\in Sp(V)$. More precisely, $B\in G^\d_{<0}$. We have $B\ne1$ since $\sum_kx_k\m(e_{-2},x_2)e_0\ne0$. Hence
$B\n G_{\ge0}^\d$. We have
$$\align&BA(\sum_kx_k)-AB(\sum_kx'_k)\\&
=\sum_{k\ne-2,0}Ax_{k-2}+(Ax_{-4}+(e_0,Ax_{-2})e_{-2})+(Ax_{-2}+(e_{-2},Ax_0)e_0)\\&
-\sum_{k\ne-2,0}Ax_k-(Ax_{-2}+(e_0,x_0)Ae_{-2})-(Ax_0+(e_{-2},x_2)Ae_0)\\&
=-(Ae_0,x_{-2})e_{-2}-(Ae_{-2},x_0)e_0-(e_0,x_0)Ae_{-2}-(e_{-2},x_2)Ae_0=0.\endalign$$
Thus $AB=BA$. We see that $A\n\fg^{\d!}_2$.

We now assume that $A:V_{-i}@>>>V_{-i+2}$ is injective for any $i\ge2$ and that for some even $n>0$, 
$A^n:V_{-n}@>>>V_n$ is not an isomorphism. The kernel of this map is the radical of the symplectic form 
$x,x'\m(x,A^nx')$ on $V_{-n}$ hence it has even codimension in $V_{-n}$; but then it also has even dimension since
$\dim V_{-n}$ is even; since this kernel is non-zero it has dimension $\ge2$. Thus we can find $e_{-n},f_{-n}$
linearly independent in $V_{-n}$ such that $A^ne_{-n}=0$, $A^nf_{-n}=0$. For $j\ge0$ we set 
$e_{2j-n}=A^je_{-n}$, $f_{2j-n}=A^jf_{-n}$. We have $e_n=0, f_n=0$. Also $e_m,f_m$ are linearly independent in 
$V_m$ if $m\le0$ is even. For $j\in[0,n]$ we have $(e_{2j-n},e_{n-2j})=0$, $(f_{2j-n},f_{n-2j})=0$,
$(e_{2j-n},f_{n-2j})=0$, $(f_{2j-n},e_{n-2j})=0$. (The last of these equalities is equivalent to
$(A^jf_{-n},A^{n-j}f_{-n})=0$ that is to $(f_{-n},A^nf_{-n})=0$ which follows from $A^nf_{-n}=0$. The other three
equalities are proved in a similar way.) Define $B\in\End(V)$ by 
$$\align&B(\sum_kx_k)=\sum_{k\n\{2h-n;h\in[0,n-1]\}}x_k\\&+
\sum_{j\in[0,n-1]}(x_{2j-n}+(-1)^j(f_{n-2j-2},x_{2j-n+2})e_{2j-n}\\&-(-1)^j(e_{n-2j-2},x_{2j-n+2})f_{2j-n}).
\endalign$$
We have 
$$\align&(B(\sum_kx_k),B(\sum_kx'_k))-(\sum_kx_k,\sum_kx'_k)\\&=
\sum_{j\in[1,n-1]}(x_{2j-n}+(-1)^j(f_{n-2j-2},x_{2j-n+2})e_{2j-n}\\&-(-1)^j(e_{n-2j-2},x_{2j-n+2})f_{2j-n},
x'_{n-2j}+(-1)^{n-j}(f_{2j-n-2},x'_{n-2j+2})e_{n-2j}\\&-(-1)^{n-j}(e_{2j-n-2},x'_{n-2j+2})f_{n-2j})\\&
+(x_{-n}+(f_{n-2},x_{2-n})e_{-n}-(e_{n-2},x_{2-n})f_{-n},x'_n)\\&
+(x_n,x'_{-n}+(f_{n-2},x'_{2-n})e_{-n}-(e_{n-2},x'_{2-n})f_{-n})\endalign$$
$$\align&+\sum_{k\ne-n,2-n,...,n-2,n}(x_k,x'_{-k})-\sum_k(x_k,x'_{-k})\\&
=\sum_{j\in[1,n]}(-1)^{n-j}(x_{2j-n},(f_{2j-n-2},x'_{n-2j+2})e_{n-2j}\\&-(e_{2j-n-2},x'_{n-2j+2})f_{n-2j})
+\sum_{j\in[0,n-1]}(-1)^j(f_{n-2j-2},x_{2j-n+2})e_{2j-n}\\&-(e_{n-2j-2},x_{2j-n+2})f_{2j-n},x'_{n-2j})\\&
=\sum_{j\in[1,n]}(-1)^{n-j}(x_{2j-n},e_{n-2j})(f_{2j-n-2},x'_{n-2j+2})\\&
-\sum_{j\in[1,n]}(-1)^{n-j}(x_{2j-n},f_{n-2j})(e_{2j-n-2},x'_{n-2j+2})\\&
+\sum_{j\in[0,n-1]}(-1)^j(f_{n-2j-2},x_{2j-n+2})(e_{2j-n},x'_{n-2j})\\&
-\sum_{j\in[0,n-1]}(-1)^j(e_{n-2j-2},x_{2j-n+2})(f_{2j-n},x'_{n-2j})\endalign$$
$$\align&=-\sum_{j\in[1,n]}(-1)^j(e_{n-2j},x_{2j-n})(f_{2j-n-2},x'_{n-2j+2})\\&
+\sum_{j\in[1,n]}(-1)^j(e_{2j-n-2},x'_{n-2j+2})(f_{n-2j},x_{2j-n})\\&
+\sum_{j\in[1,n]}(-1)^{j-1}(e_{2j-n-2},x'_{n-2j+2})(f_{n-2j},x_{2j-n})\\&
-\sum_{j\in[1,n]}(-1)^{j-1}(e_{n-2j},x_{2j-n})(f_{2j-n-2},x'_{n-2j+2})=0.\endalign$$
Thus $B\in Sp(V)$. More precisely, $B\in G^\d_{<0}$. We have $B\ne1$ since 
$x_1\m(f_{-1},x_1)e_{-1}-(e_{-1},x_1))f_{-1})$ is $\ne0$. Hence $B\n G_{\ge0}^\d$. We have
$$\align&BA(\sum_kx_k)-AB(\sum_kx_k)\\&
=\sum_{k\n\{2h-n;h\in[0,n-1]\}}Ax_{k-2}+
\sum_{j\in[0,n-1]}(Ax_{-n-2+2j}\\&+(-1)^j(f_{n-2j-2},Ax_{2j-n})e_{2j-n}
-(-1)^j(e_{n-2j-2},Ax_{2j-n})f_{2j-n})\\&
-\sum_{k\n\{2h-n;h\in[0,n-1]\}}Ax_k\\&-\sum_{j\in[0,n-1]}
(Ax_{2j-n}+(-1)^j(f_{n-2j-2},x_{2j-n+2})Ae_{2j-n}\\&-(-1)^j(e_{n-2j-2},x_{2j-n+2})Af_{2j-n})\endalign$$
$$\align&=\sum_{j\in[0,n-1]}(-1)^j((f_{n-2j-2},Ax_{2j-n})e_{2j-n}-(e_{n-2j-2},Ax_{2j-n})f_{2j-n})\\&
-\sum_{j\in[0,n-1]}(-1)^j((f_{n-2j-2},x_{2j-n+2})Ae_{2j-n}-(e_{n-2j-2},x_{2j-n+2})Af_{2j-n})\\&
=-\sum_{j\in[0,n-1]}(-1)^{j-1}((f_{n-2j},x_{2j-n})e_{2j-n}+(e_{n-2j},x_{2j-n})f_{2j-n})\\&
+\sum_{j\in[0,n-1]}(-1)^j(-(f_{n-2j-2},x_{2j-n+2})e_{2j-n+2}+(e_{n-2j-2},x_{2j-n+2})f_{2j-n+2})\\&
=\sum_{j\in[0,n-1]}(-1)^{j-1}((f_{n-2j},x_{2j-n})e_{2j-n}-(e_{n-2j},x_{2j-n})f_{2j-n})\\&
+\sum_{j\in[1,n]}(-1)^{j-1}(-(f_{n-2j},x_{2j-n})e_{2j-n}+(e_{n-2j},x_{2j-n})f_{2j-n})\\&
=-(f_n,x_{-n})e_{-n}-(e_n,x_{-n})f_{-n}+(-1)^n(f_{-n},x_n)e_n-(-1)^n(e_{-n},x_n)f_n=0\endalign$$
since $f_n=0, e_n=0$. Thus $AB=BA$. We see that $A\n\fg^{\d!}_2$.

We now assume that $A:V_{-i}@>>>V_{-i+2}$ is injective for any $i\ge2$ and that for some odd $n>0$, 
$A^n:V_{-n}@>>>V_n$ is not an isomorphism. We can find $e_{-n}\in V_{-n}-\{0\}$ such that $A^ne_{-n}=0$. For any 
$j\ge0$ we set $e_{2j-n}=A^je_{-n}\in V_{2j-n}$. Note that $e_n=0$. Also $e_m\ne0$ if $m\le0$ is odd. We have
$(e_{2j-n},e_{n-2j})=0$ for $j\in[0,n]$. Indeed we must show that $(A^je_{-n},A^{n-j}e_{-n})=0$ that is,
$(e_{-n},A^ne_{-n})=0$. This follows from $A^ne_{-n}=0$. Define $B\in\End(V)$ by   
$$\align&B(\sum_kx_k)\\&=\sum_{k\n\{2h-n;h\in[0,n-1]\}}x_k
+\sum_{j\in[0,n-1]}(x_{2j-n}+(-1)^j(e_{n-2j-2},x_{2j-n+2})e_{2j-n}).\endalign$$
We have 
$$\align&(B(\sum_kx_k),B(\sum_kx'_k))-(\sum_kx_k,\sum_kx'_k)\\&=
\sum_{j\in[1,n-1]}(x_{2j-n}+(-1)^j(e_{n-2j-2},x_{2j-n+2})e_{2j-n},\\&x'_{n-2j}+
(-1)^{n-j}(e_{2j-n-2},x'_{n-2j+2})e_{n-2j})+
(x_{-n}+(e_{n-2},x_{2-n})e_{-n},x'_n)\\&+(x_n,x'_{-n}+(e_{n-2},x'_{2-n})e_{-n})
+\sum_{k\ne0,1,...,n-1}(x_{-k},x'_k)-\sum_k(x_{-k},x'_k)\endalign$$
$$\align&=\sum_{j\in[1,n]}(-1)^{n-j}(x_{2j-n},e_{n-2j})(e_{2j-n-2},x'_{n-2j+2})\\&
+\sum_{j\in[0,n-1]}(-1)^j((e_{n-2j-2},x_{2j-n+2})(e_{2j-n},x'_{n-2j})\\&
=\sum_{j\in[1,n]}(-1)^j(e_{n-2j},x_{2j-n}))(e_{2j-n-2},x'_{n-2j+2})\\&
+\sum_{j\in[1,n]}(-1)^{j-1}((e_{n-2j},x_{2j-n})(e_{2j-n-2},x'_{n-2j+2})=0.\endalign$$
Thus $B\in Sp(V)$. More precisely, $B\in G^\d_{<0}$. We have $B\ne1$ since $x_1\m\pm(e_{-1},x_1)e_{-1}$ is $\ne0$. 
Hence $B\n G_{\ge0}^\d$. We have 
$$\align&BA(\sum_kx_k)-AB(\sum_kx'_k)=\sum_{k\n\{2h-n;h\in[0,n-1]\}}Ax_{k-2}\\&
+\sum_{j\in[0,n-1]}(Ax_{2j-n-2}+(-1)^j(e_{n-2j-2},Ax_{2j-n})e_{2j-n})\\&
-\sum_{k\n\{2h-n;h\in[0,n-1]\}}Ax_k\\&
-\sum_{j\in[0,n-1]}(Ax_{2j-n}+(-1)^j(e_{n-2j-2},x_{2j-n+2})Ae_{2j-n})\\&
=-\sum_{j\in[0,n-1]}(-1)^j(Ae_{n-2j-2},x_{2j-n})e_{2j-n}\\&
-\sum_{j\in[0,n-1]}(-1)^j(e_{n-2j-2},x_{2j-n+2})Ae_{2j-n}\\&
=-\sum_{j\in[0,n-1]}(-1)^j(e_{n-2j},x_{2j-n})e_{2j-n}\\&
+\sum_{j\in[1,n]}(-1)^j(e_{n-2j},x_{2j-n})e_{2j-n}\\&
=-(e_n,x_{-n})e_{-n}+(-1)^n(e_{-n},x_n)e_n=0.\endalign$$
Thus, $AB=BA$. We see that $A\n\fg^{\d!}_2$.

We have shown that $\fs(V)_2-\fs(V)_2^0\sub\fg^\d_2-\fg^{\d!}_2$.

Conversely, let $A\in\fs(V)_2^0$. Let $B\in G$ be such that $AB=BA$. It is enough to show that 
$B\in G_{\ge0}^\d$. Since $A\in\End(V)_2^0$ we see from the proof in 1.6 that $B(V_{\ge i})=V_{\ge i}$ for 
any $i\in\ZZ$. In particular, $B\in G_{\ge0}^\d$. This completes the proof.

\subhead 1.8\endsubhead
We prove 1.5(a). Generally, $x_k$ will denote an element of $V_k$. 
Let $A\in\fo(V)_2-\fo(V)_2^0$.

Assume that $A:V_{-i}@>>>V_{-i+2}$ is not injective for some $i>2$. Then $A:V_{i-2}@>>>V_i$ is not surjective and
since $\dim V_{i-2}\ge\dim V_i$, we see that $A:V_{i-2}@>>>V_i$ is not injective. We can find 
$e_{-i}\in V_{-i}-\{0\}$ such that $Ae_{-i}=0$. We can find $e_{i-2}\in V_{i-2}-\{0\}$ such that $Ae_{i-2}=0$.
Define $B\in\End(V)$ by
$$B(\sum_kx_k)=\sum_{k\ne-i,i-2}x_k+(x_{-i}+(e_{i-2},x_{-i+2})e_{-i})+(x_{i-2}-(e_{-i},x_i)e_{i-2}).$$
We have
$$\align&QB(\sum_kx_k)-Q(\sum_kx_k)
=Q(x_0)+(x_{-i}+(e_{i-2},x_{-i+2})e_{-i},x_i)\\&+(x_{i-2}-(e_{-i},x_i)e_{i-2},x_{-i+2})
+\sum_{k>0,k\ne i-2,i}(x_{-k},x_k)-Q(x_0)-\sum_{k>0}(x_{-k},x_k)\\&
=((e_{i-2},x_{-i+2})e_{-i},x_i)-((e_{-i},x_i)e_{i-2},x_{-i+2})\\&
=(e_{i-2},x_{-i+2})(e_{-i},x_i)-(e_{-i},x_i)(e_{i-2},x_{-i+2})=0.\endalign$$
Hence $B\in O(V)$. More precisely, $B\in G^\d_{<0}$. We have $B\ne1$ since $x_{-i+2}\m(e_{i-2},x_{-i+2})e_{-i}$ is
$\ne0$. Hence $B\n G_{\ge0}^\d$. We have 
$$\align&BA(\sum_kx_k)-AB(\sum_kx_k)\\&
=\sum_{k\ne-i,i-2}Ax_{k-2}+(Ax_{-i-2}+(e_{i-2},Ax_{-i})e_{-i})+(Ax_{i-4}-(e_{-i},Ax_{i-2})e_{i-2})\\&
-\sum_{k\ne-i,i-2}Ax_k-(Ax_{-i}+(e_{i-2},x_{-i+2})Ae_{-i})-(Ax_{i-2}-(e_{-i},x_i)Ae_{i-2})\\&
=(e_{i-2},Ax_{-i})e_{-i}-(e_{-i},Ax_{i-2})e_{i-2})-(e_{i-2},x_{-i+2})Ae_{-i}+(e_{-i},x_i)Ae_{i-2}\\&
=-(Ae_{i-2},x_{-i})e_{-i}+(Ae_{-i},x_{i-2})e_{i-2})-(e_{i-2},x_{-i+2})Ae_{-i}+(e_{-i},x_i)Ae_{i-2}
=0.\endalign$$
Hence $AB=BA$. We see that $A\n\fg^{\d!}_2$.

Next we assume that $A:V_{-2}@>>>V_0$ is not injective that is,\lb $K:=\ker(A:V_{-2}@>>>V_0)\ne0$. We can find 
$e_{-2}\in V_{-2}-\{0\}$ such that $Ae_{-2}=0$. Note that $K':=\ker(A:V_0@>>>V_2)$ is $\ne0$. (If $\dim V$ is
odd and $p=2$ we have $0\ne\car\sub K'$. If $\dim V$ is even or if $p\ne2$ then $A:V_0@>>>V_2$ is the transpose of
$-A:V_{-2}@>>>V_0$ hence is not surjective. But $\dim V_0\ge\dim V_2$ hence again $K'\ne0$.) We can find 
$e_0\in V_0-\{0\}$ such that $Ae_0=0$. Let $I'=\{x\in V_2;(x,K)=0\}$. Note that $I'\ne V_2$. Define a linear 
function $f:V_2@>>>\kk$ by $f(x_2)=(e_{-2},x_2)\sqrt{Q(e_0)}$ where $\sqrt{Q(e_0)}$ is a fixed square root of 
$Q(e_0)$. Now $f=0$ on $I'$ since $e_{-2}\in K$. Hence $f$ induces a linear function $f':V_2/I'@>>>\kk$. Note 
that $K=(V_2/I')^*$ canonically. There is a unique linear function $\g':V_2/I'@>>>K$ such that 
$(\g'(x_2),x'_2)=-f'(x_2)f'(x'_2)$ for all $x_2,x'_2\in V_2/I'$. Hence there is a unique linear function 
$\g:V_2@>>>V_{-2}$ such that $(\g(x_2),x'_2)=-f(x_2)f(x'_2)$ for all $x_2,x'_2$ in $V_2$ and $\g(V_2)\sub K$, 
$\g(I')=0$. Hence $A\g(V_2)=0$. Since $AV_0\sub I'$, we have $\g AV_0=0$. Define $B\in\End(V)$ by 
$$B(\sum_kx_k)=\sum_{k\ne-2,0}x_k+(x_{-2}+(e_0,x_0)e_{-2}+\g(x_2))+(x_0-(e_{-2},x_2)e_0).$$
We have 
$$\align&QB(\sum_kx_k)-Q(\sum_kx_k)\\&=
Q(x_0)+(e_{-2},x_2)^2Q(e_0)-(e_0,x_0)(e_{-2},x_2)+(x_{-2},x_2)\\&+(e_0,x_0)(e_{-2},x_2)+(\g(x_2),x_2)
+\sum_{k>0,k\ne2}(x_{-k},x_k)-Q(x_0)-\sum_{k>0}(x_{-k},x_k)\\&=f(x_2)^2+(\g(x_2),x_2)=0.\endalign$$
Thus $B\in O(V)$. More precisely, $B\in G^\d_{<0}$. We have $B\ne1$ since $x_2\m(e_{-2},x_2)e_0$ is $\ne0$. 
Hence $B\n G_{\ge0}^\d$. We have 
$$\align&BA(\sum_kx_k)-AB(\sum_kx_k)\\&=
\sum_{k\ne-2,0}Ax_{k-2}+(Ax_{-4}+(e_0,Ax_{-2})e_{-2}+\g(Ax_0))+(Ax_{-2}-(e_{-2},Ax_0)e_0)\\&
-\sum_{k\ne-2,0}Ax_k-(Ax_{-2}+(e_0,x_0)Ae_{-2}+A\g(x_2))-(Ax_0-(e_{-2},x_2)Ae_0)\\&
=-(Ae_0,x_{-2},e_0)e_{-2}+\g(Ax_0)-(Ae_{-2},x_0)e_0-(e_0,x_0)Ae_{-2}-A\g(x_2)\\&+(e_{-2},x_2)Ae_0
=\g(Ax_0)-A\g(x_2)=0.\endalign$$
Thus $AB=BA$. We see that $A\n\fg^{\d!}_2$.

We now assume that $A:V_{-i}@>>>V_{-i+2}$ is injective for any $i\ge2$ and that for some $n>0$, 
and some $\x\in A^n(V_{-2n})-\{0\}$ we have $(\x,A^n(V_{-2n})=0$, $Q(\x)=0$. We can 
write $\x=A^ne_{-2n}$ for a unique $e_{-2n}\in V_{-2n}-\{0\}$. For any $j\ge0$ we set 
$e_{-2n+2j}=A^je_{-2n}\in V_{-2n+2j}$. Thus $e_0=\x$ and $(e_0,A^nV_{-2n})=0$, $Q(e_0)=0$. We show that 
$e_{2n}=0$. Indeed, $(V_{-2n},e_{2n})=(V_{-2n},A^ne_0)=\pm(A^nV_{-2n},e_0)=0$. 
For $j\in[0,2n]$ we have:

$(e_{-2n+2j},e_{2n-2j})=0$.
\nl
This follows from $(A^je_{-2n},A^{2n-j}e_{-2n})=\pm(A^{2n}e_{-2n},e_{-2n})=(e_{2n},e_{-2n})=0$.

Define $B\in\End(V)$ by
$$\align&B(\sum_kx_k)=\sum_{k\n\{-2n+2h;h\in[0,2n-1]\}}x_k\\&+
\sum_{j\in[0,2n-1]}(x_{-2n+2j}+(-1)^j(e_{2n-2j-2},x_{-2n+2j+2})e_{-2n+2j}).\endalign$$
We have 
$$\align&QB(\sum_kx_k)-Q(\sum_kx_k)\\&=
Q(x_0+(e_{-2},x_2)e_0)+\sum_{j\in[1,n-1]}
(x_{-2n+2j}+(-1)^j(e_{2n-2j-2},x_{-2n+2j+2})e_{-2n+2j},\\&
x_{2n-2j}+(-1)^{2n-j}(e_{-2n+2j-2},x_{2n-2j+2})e_{2n-2j})\\&
+(x_{-2n}+(e_{2n-2},x_{-2n+2})e_{-2n},x_{2n})\\&
+\sum_{k>0;k\ne2,4,...,2n}(x_{-k},x_k)-Q(x_0)-\sum_{k>0}(x_{-k},x_k)\endalign$$
$$\align&=Q(e_0)+(e_0,x_0)(e_{-2},x_2)\\&
+\sum_{j\in[1,n-1]}(-1)^j(e_{2n-2j},x_{-2n+2j})(e_{-2n+2j-2},x_{2n-2j+2})\\&
+\sum_{j\in[1,n-1]}(-1)^j(e_{2n-2j-2},x_{-2n+2j+2})(e_{-2n+2j},x_{2n-2j})\\&
+(e_{2n-2},x_{-2n+2})(e_{-2n},x_{2n})-Q(x_0)-\sum_{k>0}(x_{-k},x_k)\endalign$$
$$\align&=\sum_{j\in[2,n-1]}(-1)^j(e_{2n-2j},x_{-2n+2j})(e_{-2n+2j-2},x_{2n-2j+2})\\&
+\sum_{j\in[1,n-2]}(-1)^j(e_{2n-2j-2},x_{-2n+2j+2})(e_{-2n+2j},x_{2n-2j})\\&
=\sum_{j\in[2,n-1]}(-1)^j(e_{2n-2j},x_{-2n+2j})(e_{-2n+2j-2},x_{2n-2j+2})\\&
+\sum_{j\in[2,n-1]}(-1)^{j-1}(e_{2n-2j},x_{-2n+2j})(e_{-2n+2j-2},x_{2n-2j+2})=0.\endalign$$
Hence $B\in SO(V)$. More precisely, $B\in G^\d_{<0}$. We have $B\ne1$ since $x_2\m(e_{-2},x_2)e_0$ is $\ne0$.
Hence $B\n G_{\ge0}^\d$. We have  
$$\align&BA(\sum_kx_k)-AB(\sum_kx_k)\\&=
\sum_{k\n\{-2n+2h;h\in[0,2n-1]\}}Ax_{k-2}\\&+
\sum_{j\in[0,2n-1]}(Ax_{-2n+2j-2}+(-1)^j(Ax_{-2n+2j},e_{2n-2j-2})e_{-2n+2j})\\&
-\sum_{k\n\{-2n+2h;h\in[0,2n-1]\}}Ax_k\\&
-\sum_{j\in[0,2n-1]}(Ax_{-2n+2j}+(-1)^j(x_{-2n+2j+2},e_{2n-2j-2})e_{-2n+2j+2})\endalign$$
$$\align&=\sum_{k}Ax_k+\sum_{j\in[0,2n-1]}(-1)^j(Ax_{-2n+2j},e_{2n-2j-2})e_{-2n+2j})\\&
-\sum_{k}Ax_k-\sum_{j\in[0,2n-1]}(-1)^j(x_{-2n+2j+2},e_{2n-2j-2})e_{-2n+2j+2}\\&
=-\sum_{j\in[0,2n-1]}(-1)^j(x_{-2n+2j},Ae_{2n-2j-2})e_{-2n+2j})\\&
-\sum_{j\in[1,2n]}(-1)^{j-1}(x_{-2n+2j},e_{2n-2j})e_{-2n+2j}\\&
=-(x_{-2n},e_{2n})e_{-2n}-(-1)^n(x_{2n},e_{-2n})e_{2n}=0 \endalign$$
since $e_{2n}=0$.
Hence $AB=BA$. We see that $A\n\fg^{\d!}_2$.

We now assume that $A:V_{-i}@>>>V_{-i+2}$ is injective for any $i\ge2$ and that for some odd $n>0$,
$A^n:V_{-n}@>>>V_n$ is not an isomorphism. The kernel of this map is the radical of the symplectic form 
$x,x'\m(x,A^nx')$ on $V_{-n}$ hence it has even codimension in $V_{-n}$; but then it also has even dimension since
$\dim V_{-n}$ is even; since this kernel is non-zero it has dimension $\ge2$. Thus we can find $e_{-n},f_{-n}$
linearly independent in $V_{-n}$ such that $A^ne_{-n}=0$, $A^nf_{-n}=0$. For $j\ge0$ we set 
$e_{2j-n}=A^je_{-n}$, $f_{2j-n}=A^jf_{-n}$. We have $e_n=0,f_n=0$. Also $e_m,f_m$ are linearly independent in 
$V_m$ if $m<0$ is odd. For $j\in[0,n]$ we have $(e_{2j-n},e_{n-2j})=0$, $(f_{2j-n},f_{n-2j})=0$,
$(e_{2j-n},f_{n-2j})=0$, $(f_{2j-n},e_{n-2j})=0$. (The last of these equalities is equivalent to
$(A^jf_{-n},A^{n-j}f_{-n})=0$ that is, to $(f_{-n},A^nf_{-n})=0$ which follows from $A^nf_{-n}=0$. The other three
equalities are proved in a similar way.) Define $B\in\End(V)$ by   
$$\align&B(\sum_kx_k)=\sum_{k\n\{2h-n;h\in[0,n-1]\}}x_k+\\&
\sum_{j\in[0,n-1]}(x_{2j-n}+(-1)^j(f_{n-2j-2},x_{2j-n+2})e_{2j-n}\\&-(-1)^j(e_{n-2j-2},x_{2j-n+2})f_{2j-n}).
\endalign$$
We have   
$$\align&QB(\sum_kx_k)-Q(\sum_kx_k)\\&
=Q(x_0)+\sum_{j\in[1,(n-1)/2]}(x_{2j-n}+(-1)^j(f_{n-2j-2},x_{2j-n+2})e_{2j-n}\\&-
(-1)^j(e_{n-2j-2},x_{2j-n+2})f_{2j-n},\\&
x_{n-2j}+(-1)^{n-j}(f_{2j-n-2},x_{n-2j+2})e_{n-2j}-(-1)^{n-j}(e_{2j-n-2},x_{n-2j+2})f_{n-2j})\\&
+(x_{-n}+(f_{n-2},x_{2-n})e_{-n}-(e_{n-2},x_{2-n})f_{-n},x_n)\\&+
\sum_{k>0;k\ne-n,2-n,...,n-2}(x_{-k},x_k)+Q(x_0)-\sum_{k>0}(x_{-k},x_k)-Q(x_0)\endalign$$
$$\align&=\sum_{j\in[1,(n-1)/2]}(x_{2j-n},(-1)^{n-j}(f_{2j-n-2},x_{n-2j+2})e_{n-2j}\\&-(-1)^{n-j}
(e_{2j-n-2},x_{n-2j+2})f_{n-2j})\\&+\sum_{j\in[0,(n-1)/2]}((-1)^j(f_{n-2j-2},x_{2j-n+2})e_{2j-n}\\&
-(-1)^j(e_{n-2j-2},x_{2j-n+2})f_{2j-n},x_{n-2j})\\&
=\sum_{j\in[1,(n-1)/2]}(-1)^{n-j}(e_{n-2j},x_{2j-n})(f_{2j-n-2},x_{n-2j+2})\\&
-\sum_{j\in[1,(n-1)/2]}(-1)^{n-j}(f_{n-2j},x_{2j-n})(e_{2j-n-2},x_{n-2j+2})\\&
-\sum_{j\in[0,(n-1)/2]}(-1)^j(f_{2j-n},x_{n-2j})(e_{n-2j-2},x_{2j-n+2})\\&
+\sum_{j\in[0,(n-1)/2]}(-1)^j(e_{2j-n},x_{n-2j})(f_{n-2j-2},x_{2j-n+2})\endalign$$
$$\align&=\sum_{j\in[1,(n-1)/2]}(-1)^{n-j}(e_{n-2j},x_{2j-n})(f_{2j-n-2},x_{n-2j+2})\\&
-\sum_{j\in[1,(n-1)/2]}(-1)^{n-j}(f_{n-2j},x_{2j-n})(e_{2j-n-2},x_{n-2j+2})\\&
-\sum_{j\in[1,(n+1)/2]}(-1)^{j-1}(f_{2j-n-2},x_{n+2-2j})(e_{n-2j},x_{2j-n})\\&
+\sum_{j\in[1,(n+1)/2]}(-1)^{j-1}(e_{2j-n-2},x_{n+2-2j})(f_{n-2j},x_{2j-n})\\&
=(-1)^{(n-1)/2}(e_{-1},x_1)(f_{-1},x_1)-(-1)^{(n-1)/2}(f_{-1},x_1)(e_{-1},x_1)=0.\endalign$$
Thus, $B\in O(V)$. More precisely, $B\in G^\d_{<0}$. We have $B\ne1$ since 
$x_1\m\pm((f_{-1},x_1)e_{-1}-(e_{-1},x_1))f_{-1}))$ is $\ne0$. Hence $B\n G_{\ge0}^\d$. We have 
$$\align&BA(\sum_kx_k)-AB(\sum_kx_k)=\sum_{k\n\{2h-n;h\in[0,n-1]\}}Ax_{k-2}\\&+
\sum_{j\in[0,n-1]}(Ax_{-n-2+2j}+(-1)^j(f_{n-2j-2},Ax_{2j-n})e_{2j-n}-\\&(-1)^j(e_{n-2j-2},Ax_{2j-n})
f_{2j-n})-\sum_{k\n\{2h-n;h\in[0,n-1]\}}Ax_k\\&
-\sum_{j\in[0,n-1]}(Ax_{2j-n}+(-1)^j(f_{n-2j-2},x_{2j-n+2})Ae_{2j-n}\\&
-(-1)^j(e_{n-2j-2},x_{2j-n+2})Af_{2j-n})\endalign$$
$$\align&=\sum_{j\in[0,n-1]}(-1)^j((f_{n-2j-2},Ax_{2j-n})e_{2j-n}-(e_{n-2j-2},Ax_{2j-nj})f_{2j-n})\\&
-\sum_{j\in[0,n-1]}(-1)^j((f_{n-2j-2},x_{2j-n+2})Ae_{2j-n}-(e_{n-2j-2},x_{2j-n+2})Af_{2j-n})\\&
=-\sum_{j\in[0,n-1]}(-1)^{j-1}((f_{n-2j},x_{2j-n})e_{2j-n}+(e_{n-2j},x_{2j-n})f_{2j-n})\\&
+\sum_{j\in[0,n-1]}(-1)^j(-(f_{n-2j-2},x_{2j-n+2})e_{2j-n+2}+(e_{n-2j-2},x_{2j-n+2})f_{2j-n+2})\\&
=\sum_{j\in[0,n-1]}(-1)^{j-1}((f_{n-2j},x_{2j-n})e_{2j-n}-(e_{n-2j},x_{2j-n})f_{2j-n})\\&
+\sum_{j\in[1,n]}(-1)^{j-1}(-(f_{n-2j},x_{2j-n})e_{2j-n}+(e_{n-2j},x_{2j-n})f_{2j-n})\\&
=-(f_n,x_{-n})e_{-n}+(e_n,x_{-n})f_{-n}+(-1)^n(f_{-n},x_n)e_n-(-1)^n(e_{-n},x_n)f_n=0\endalign$$
since $f_n=0$, $e_n=0$.
Thus $AB=BA$. We see that $A\n\fg^{\d!}_2$.

We have shown that $\fo(V)-\fo(V)_2^0\sub\fg^\d_2-\fg^{\d!}_2$. 

Conversely, let $A\in\fo(V)_2^0$. Let $B\in G$ be such that $AB=BA$. It is enough to show that $B\in G_{\ge0}^\d$.
We argue by induction on $\dim V$. If $V=0$ the result is 
clear. We now assume $V\ne0$. Let $m$ be the largest integer $\ge0$ such that $V_m\ne0$. If $m=0$ we have
$G_{\ge0}^\d=G$ and the result is clear. Assume now that $m\ge1$. 
If $m$ is odd we have $A^mV=V_m$, $\ker(A^m:V@>>>V)=V_{\ge-m+1}$. Since $BA=AB$, the image and kernel of $A^m$
are $B$-stable. Hence $B(V_m)=V_m$ and
$B(V_{\ge-m+1})=V_{\ge-m+1}$. Hence $B$ induces an automorphism $B'\in SO(V')$ where 
$V'=V_{\ge-m+1}/V_m$, a vector spaces with a nondegenerate quadratic form induced by $Q$. We have 
canonically $V'=V_{-m+1}\op V_{-m+2}\op\do\op V_{m-1}$ and $\fo(V')_2,\fo(V')_2^0$ are defined in terms of this
($o$-good) grading. Now $A$ induces an element $A'\in\fo(V')_2^0$ and we
have $B'A'=A'B'$. By the induction hypothesis, for any $i\in[-m+1,m-1]$, the subspace $V_i+V_{i+1}+\do+V_{m-1}$ of
$V'$ is $B'$-stable. Hence the subspace $V_{\ge i}$ of $V$ is $B$-stable. We see that $B\in G_{\ge0}^\d$. 
Next we assume that $m$ is even. We have $V_{\ge-m+1}=\{x\in V; A^m(x)=0,Q(A^{m/2}x)=0\}$ (we use
that $A\in\fo(V)_2^0$). Since $B$ commutes with $A$ and preserves $Q$ we see that $B$ preserves the subspace
$\{x\in V;A^m(x)=0,Q(A^{m/2}x)=0\}$ hence $B(V_{\ge-m+1})=V_{\ge-m+1}$. We have 
$V_m=\{x\in V;(x,V_{\ge-m+1})=0,Q(x)=0\}$. Since $B$ preserves the subspace
$V_{\ge-m+1}$ and $B$ preserves $Q$ and $(,)$ we see that $B(V_m)=V_m$. 
Hence $B$ induces an automorphism $B'\in SO(V')$ where 
$V'=(V_{\ge-m+1})/V_m$, a vector space with a nondegenerate quadratic form induced by $Q$. We have 
canonically $V'=V_{-m+1}\op V_{-m+2}\op\do\op V_{m-1}$ and $\fo(V')_2,\fo(V')_2^0$ are defined in terms of this
($o$-good) grading. Now $A$ induces an element $A'\in\fo(V')_2^0$ and we
have $B'A'=A'B'$. By the induction hypothesis, for any $i\in[-m+1,m-1]$, the subspace $V_i+V_{i+1}+\do+V_{m-1}$ of
$V'$ is $B'$-stable. Hence the subspace $V_{\ge i}$ of $V$ is $B$-stable. We see that
$B\in G_{\ge0}^\d$. This completes the proof.

\subhead 1.9\endsubhead
Let $\d\in\fD_G$. We describe the set of $G^\d_0$-orbits of $\fg_2^{\d!}$ (see 1.2) in the cases considered in
1.3-1.5. If $G,\fg$ are as in 1.3 or as in 1.4 (with $p\ne2$) or 1.5 (with $p\ne2$) then, using 1.3(a), 1.4(a),
1.5(a), we see that $G^\d_0$ acts transitively on $\fg_2^{\d!}$. If $V,G,\fg$ are as in 1.4 (with $p=2$) then the
set of $G^\d_0$-orbits of $\fg_2^{\d!}=\fs(V)_2^0$ is a set (with cardinal a power of $2$) described in
\cite{\LUU, p.478}. In the rest of this subsection we assume that $V,Q,(,),G,\fg$ are as in 1.5 and $p=2$. Let 
$(V_i)$ be the $o$-good grading of $V$ corresponding to $\d$. For any $n\ge0$ let $d_n=\dim V_{-2n}$. Let 
$M=\{n\ge0;d_n=\text{odd}\}$. If $M=\em$ then $G^\d_0$ acts transitively on $\fg_2^{\d!}=\fo(V)_2^0$. Now assume 
that $M\ne\em$. We write the elements of $M$ in increasing order $n_0<n_1<\do<n_t$. 
Let $X$ be the set of all functions $f:\{1,2,\do,t\}@>>>\{0,1\}$ such that:

(i)  $f(i)=1$ if $n_i-n_{i-1}\ge2$;

(ii) $f(i)=0$ if $n_i-n_{i-1}=1$ and $d_{n_{i-1}}=d_{n_i}$.
\nl
Note that $|X|=2^\a$ where $\a=|\{i\in[1,t];n_{i-1}=n_i-1,d_{n_{i-1}}>d_{n_i}\}|$. 

For any $A\in\fo(V)_2^0$ and any $n\in M$ let $L_n^A$ be the radical of the restriction of $(,)$ to 
$K_n^A:=A^n(V_{-2n})$ (a line). 

Note that: if $0\le s<q<r$, $s\in M,q\n M,r\in M$ then $L^A_s\ne L^A_r$. (Indeed, let $x\in L^A_s-0$. Then 
$(x,K^A_q)=0$ hence $x\n K^A_q$ hence $x\n K^A_r$ and $x\n L^A_r$.) If $0\le s<q<r$, $s\in M,q\in M,r\in M$ and 
$L_s^A=L_r^A$ then $L^A_s=L_q^A=L_r^A$. (Indeed, let $x\in L_s^A-0$. Then $(x,K_q^A)=0$. But 
$x\in L^A_r\sub K^A_r\sub K^A_q$ hence $x\in L^A_q$. Thus $L^A_s=L_q^A$.) If $s\in M,q\in M$, $K_s^A=K_q^A$ then 
clearly $L_s^A=L_q^A$.

We define $f_A:\{1,2,\do,t\}@>>>\{0,1\}$ by 
$f_A(i)=0$ if $L^A_{n_i}=L^A_{n_{i-1}}$, $f_A(i)=1$ if $L^A_{n_i}\ne L^A_{n_{i-1}}$. From the previous paragraph 
we see that $f_A\in X$ and that for any $i,j$ in 
$[0,t]$, $f_A$ determines whether $L^A_{n_i},L^A_{n_j}$ are equal or not.

For any $f\in X$ we set ${}^f\fo(V)_2^0=\{A\in\fo(V)_2^0;f_A=f\}$. Then the subsets ${}^f\fo(V)_2^0$ ($f\in X)$ 
are exactly the orbits of $G^\d_0$ on $\fg_2^{\d!}=\fo(V)_2^0$. 

\head 2. The pieces in the unipotent variety of $G$\endhead
\subhead 2.1\endsubhead
Given $\d,\d'$ in $\fD_G$ we write $\d\si\d'$ if for any $i\in\NN$ we have $G_{\ge i}^\d=G_{\ge i}^{\d'}$ or 
equivalently $\fg_{\ge i}^\d=\fg_{\ge i}^{\d'}$. This is an equivalence relation on $\fD_G$. Let $D_G$ be the set
of equivalence classes. The conjugation $G$-action on $\fD_G$ induces a $G$-action on $D_G$.
If $\D\in D_G$, $i\in\ZZ$, we can write $G^\D_{\ge i}$, $\fg^\D_{\ge i}$ instead of 
$G^\d_{\ge i},\fg^\d_{\ge i}$ where $\d\in\D$. If $\D\in D_G$, we have an action of $G_{\ge0}^\D$ on $\D$ given by
$g:\d\m\d'$ where $\d'(a)=g\d(a)g\i$ for all $a\in\kk^*$. We show:

(a) {\it the conjugation action of $G_{\ge0}^\D$ on the set of pairs $(\d,T)$ where $\d\in\D$ and $T$ is a maximal
torus of $G_{\ge0}^\D$ containing $\d(\kk^*)$ is transitive; hence the conjugation action of $G_{\ge0}^\D$ on $\D$
is transitive.}
\nl
Let $(\d,T)$, $(\d',T')$ be two pairs as above. Since $T,T'$ are conjugate in $G_{\ge0}^\D$ we can assume that
$T'=T$. It is enough to show that in this case we have $\d=\d'$. For any root $\a:T@>>>\kk^*$ of $G$ with respect
to $T$ we set 
$$\fg_\a=\{x\in\fg;\Ad(t)x=\a(t)x\qua\frl t\in T\}.$$
For any root $\a$ define $\la\d,\a\ra\in\ZZ$, $\la\d',\a\ra\in\ZZ$ by $\a(\d(a))=a^{\la\d,\a\ra}$, 
$\a(\d'(a))=a^{\la\d',\a\ra}$ for all $a\in\kk^*$. For $i>0$ we have
$$\op_{\a;\la\d,\a\ra=i}\fg_\a=\fg^\d_{\ge i}=\fg^{\d'}_{\ge i}=\op_{\a;\la\d',\a\ra=i}\fg_\a.$$
Since each $\fg_\a$ is $1$-dimensional and the sum $\sum_\a\fg_\a$ is direct it follows that 
$$\{\a;\la\d,\a\ra=i\}=\{\a;\la\d',\a\ra=i\}$$
for any $i>0$. But then we automatically have $\{\a;\la\d,\a\ra=i\}=\{\a;\la\d',\a\ra=i\}$ for any $i\in\ZZ$. It 
follows that $\la\d,\a\ra=\la\d',\a\ra$ for any root $\a$. Define $\mu\in\Hom(\kk^*,T)$ by $\mu(a)=\d'(a)\d(a)\i$
for any $a\in\kk^*$. Then $\la\mu,a\ra=0$ for any root $\a$. Thus $\mu(\kk^*)$ is contained in the centre of $G$.
Since $\d(\kk^*),\d'(\kk^*)$ are contained in $G^{der}$, we have $\mu(\kk^*)\sub G^{der}$. Thus $\mu(\kk^*)$ is 
contained in the centre of $G^{der}$, a finite group. Since $\mu(\kk^*)$ is connected we have $\mu(\kk^*)=\{1\}$ 
hence $\d'=\d$. This proves (a).

From (a) we see that 

(b) {\it the obvious map $\fD_G@>>>D_G$ induces a bijection $G\bsl\fD_G@>>>G\bsl D_G$ on the sets of 
$G$-orbits.}
\nl
In the remainder of this subsection we assume that 

(c) {\it $p>1$, $\kk$ is an algebraic closure of the field $\FF_p$ with $p$ elements
and we are given a split $\FF_p$-rational structure on $G$ with Frobenius map $F:G@>>>G$.}
\nl
This induces a split $\FF_p$-rational structure on $\fg$.
For any $\d\in\Hom(\kk^*,G)$ we define ${}^F\d:\kk^*@>>>G$ by $({}^F\d)(z)=F(\d(z^{p\i}))$ for $z\in\kk^*$. 
We show that for some $g\in G$ we have ${}^F\d(z)=g\d(z)g\i$
for all $z$. Let $T$ be a maximal torus of $G$ such that $F(t)=t^p$ for all $t\in T$.
We can find $g_1\in G$ such that $g_1\d(\kk^*)g_1\i\sub T$.
Then for $z\in\kk^*$ we have $g_1\d(z)g_1\i\sub T$ hence $F(g_1\d(z)g_1\i)=(g_1\d(z)g_1\i)^p=g_1\d(z^p)g_1\i
=F(g_1){}^F\d(z^p)F(g_1\i)$ hence ${}^F\d(z)=g\d(z)g\i$ where $g=F(g_1)\i g_1$, as claimed. 
In particular, if $\d\in\fD_G$ then ${}^F\d\in\fD_G$ and ${}^F\d$ is in the same $G$-orbit as $\d$.
From the definitions we see that if $\d\in\fD_G$ and $i\in\NN$ then $\fg_{\ge i}^{({}^F\d)}=F(\fg_{\ge i}^\d)$ 
and $G_{\ge i}^{({}^F\d)}=F(G_{\ge i}^\d)$.
In particular if $\d,\d'$ in $\fD_G$ satisfy $\d\si\d'$ then ${}^F\d\si{}^F\d'$. Thus the permutation
$\d\m{}^F\d$ of $\fD_G$ induces a permutation $\D\m{}^F\D$ of $D_G$.
This permutation maps each $G$-orbit in $D_G$ into itself.

\subhead 2.2\endsubhead
Let $\D\in D_G$, $i>0$. Let $(\d,T)$ be such that $\d\in\D$ and $T$ is a maximal torus of $G^\D_{\ge0}$ containing
$\d(\kk^*)$. We will define an isomorphism of algebraic groups
$$\Ph_{\d,T}:\fg^\D_{\ge i}/\fg^\D_{\ge i+1}@>\si>>G^\D_{\ge i}/G^\D_{\ge i+1}.$$
For any root $\a:T@>>>\kk^*$ we define the root subspaces $\fg_\a$ as in the proof of 2.1(a). Let $G_\a$ be the 
root subgroup of $G$ corresponding to $\fg_\a$. For any $\a$ we can find an isomorphism of algebraic groups 
$h_\a:\kk@>\si>>G_\a$; it is unique up to composing with multiplication by a nonzero scalar on $\kk$; by passage 
to Lie algebras, $h_\a$ gives rise to an isomorphism of vector spaces $h'_\a:\kk@>\si>>\fg_\a$. For any $j\ge0$
let $R_j$ be the set of roots $\a$ such that $\la\a,\d\ra=j$ (notation as in the proof of 2.1(a)). Then
$$\ph:\kk^{R_i}@>>>\fg^\D_{\ge i}/\fg^\D_{\ge i+1},\qua (c_\a)\m\sum_{\a\in R_i}h'_{\a}(c_\a)\mod\fg^\D_{\ge i+1},
$$
$$\ps:\kk^{R_i}@>>>G^\D_{\ge i}/G^\D_{\ge i+1},\qua(c_\a)\m\prod_{\a\in R_i}h_{\a}(c_\a)\mod G^\D_{\ge i+1}$$
are isomorphisms of algebraic groups.
(The last product is independent of the order of the factors up to $\mod G^\D_{\ge i+1}$, by the Chevalley 
commutator formula.) We set $\Ph_{\d,T}=\ps\ph\i$. Clearly $\Ph_{\d,T}$ is independent of the choice of the 
$h_\a$. Now $G_{\ge0}^\D$ acts by conjugation on $G^\D_{\ge i}/G^\D_{\ge i+1}$ and this induces an $\Ad$-action of
$G_{\ge0}^\D$ on $\fg^\D_{\ge i}/\fg^\D_{\ge i+1}$. From the definitions we see that for any $g\in G_{\ge0}^\D$ 
and any $\x\in\fg^\D_{\ge i}/\fg^\D_{\ge i+1}$ we have $\Ph_{g\d g\i,gTg\i}(\x)= g\Ph_{\d,T}(\Ad(g\i)\x)g\i$. Let
$\a\in R_i,\b\in R_j$, $j\ge0$ and let $c,c'\in\kk$. If $j>0$ we have
$h_\b(c)h_\a(c')h_\b(c)\i=h_\a(c')\mod G^\D_{\ge i+1}$, $\Ad(h_\b(c))h'_\a(c')=h'_\a(c')\mod \fg^\D_{\ge i+1}$. If
$j=0$ we have (by the Chevalley commutator formula)

$h_\b(c)h_\a(c')h_\b(c)\i=h_\a(c')\prod_{i'>0}h_{i'\b+\a}(m_{i'}c^{i'}c')\mod G^\D_{\ge i+1}$,
\nl
where $m_{i'}\in\ZZ$ are such that

$\Ad(h_\b(c))h'_\a(c')=h'_\a(c')+\sum_{i'>0}h'_{i'\b+\a}(m_{i'}c^{i'}c')$.
\nl
From these formulas we see that for any $\x\in\fg^\D_{\ge i}/\fg^\D_{\ge i+1}$ we have
$\Ph_{\d,T}(\Ad(g)\x)=g\Ph_{\d,T}(\x)g\i$ whenever $g=h_\b(c)$ with $\b\in\cup_{j\ge0}R_j$, $c\in\kk$. The same 
holds when $g\in T$ and even for any $g$ in $G_{\ge0}^\D$ since this group is generated by $T$ and by the 
$h_\b(c)$ as above. We see that for any $g\in G_{\ge0}^\D$ and any $\x\in\fg^\D_{\ge i}/\fg^\D_{\ge i+1}$ we have
$g\Ph_{\d,T}(\Ad(g\i)x)g\i=\Ph_{\d,T}(\x)$, hence $\Ph_{g\d g\i,gTg\i}(\x)=\Ph_{\d,T}(\x)$. Using this and 2.1(a)
we see that $\Ph_{\d,T}$ is independent of the choice of $\d,T$. Hence it can be denoted by $\Ph_\D$. We can 
summarize the results above as folows.

(a) {\it For any $\D\in D_G$ and any $i>0$ there is a canonical $G_{\ge0}^\D$-equivariant isomorphism of algebraic
groups $\Ph_\D:\fg^\D_{\ge i}/\fg^\D_{\ge i+1}@>\si>>G^\D_{\ge i}/G^\D_{\ge i+1}$.}

\subhead 2.3\endsubhead  
Let $\D\in D_G$. For any $\d\in\D$ the subset $\fg_2^{\d!}$ of $\fg_2^\d$ can be viewed as a subset $\Si^\d$ of
$\fg_{\ge2}^\D/\fg_{\ge3}^\D$ via the obvious isomorphism $\fg_2^\d@>\si>>\fg_{\ge2}^\D/\fg_{\ge3}^\D$. From the 
definitions, for any $g\in G^\D_{\ge0}$ we have $\Si^{g\d g\i}=\Ad(g)\Si^\d$. Here we use the $\Ad$-action of
$G^\D_{\ge0}$ on $\fg_{\ge2}^\D/\fg_{\ge3}^\D$. This action factors through an action of 
$G^\D_{\ge0}/G^\D_{\ge1}=G_0^\d$. Hence if we write $g=g_0g'$ where $g_0\in G^\d_0,g'\in G_{\ge1}^\D$ we have
$\Ad(g)(\Si^\d)=\Ad(g_0)\Si^\d=\Si^\d$ (the last equality follows from $\Ad(g_0)\fg_2^{\d!}=\fg_2^{\d!}$). Thus we
have $\Si^{g\d g\i}=\Si^\d$. Using 2.1(a) we deduce that $\Si^\d$ is independent of the choice of $\d$ in $\D$; we
will denote it by $\Si^\D$. Note that $\Si^\D$ is a subset of $\fg_{\ge2}^\D/\fg_{\ge3}^\D$ stable under the 
action of $G^\D_{\ge0}$. 

Let $\fS^\D\sub G_{\ge2}^\D$ be the inverse image of $\Si^\D$ under the composition
$G_{\ge2}^\D@>>>G_{\ge2}^\D/G_{\ge3}^\D@>\Ph_\D\i>>\fg_{\ge2}^\D/\fg_{\ge3}^\D$ where the first map is the obvious
one. Now $\fS^\D$ is stable under the conjugation action of $G^\D_{\ge0}$ on $G_{\ge2}^\D$ and $u\m u$ is a map
$$\Ps_G:\sqc_{\D\in D_G}\fS^\D@>>>\cu_G, u\m u.\tag a$$

\proclaim{Theorem 2.4} Assume that $\tG^{der}$ (see 1.1) is a product of almost simple groups of type $A,B,C,D$.
Then $\Ps_G$ is a bijection.
\endproclaim
The general case reduces easily to the case where $G$ is almost simple of type $A,B,C$ or $D$. Moreover
we can assume that $G$ is one of the groups $GL(V)$, $Sp(V)$, $SO(V)$ 
in 1.3-1.5. The proof in these cases will be given
in 2.5-2.10. We expect that the theorem holds without restriction on $G$.

We now discuss some applications of the theorem.
Let $\fA_G$ be the set of $G$-orbits on $D_G$. Using 2.1(b) and the definitions we see that $\fA_G=\fA_{G'}$ where
$G'$ is as in 1.1. In particular $\fA_G$ is a finite set which depends only on the type of $G$, not on $\kk$. For
any $\co\in\fA_G$ we consider the set

$Z_\co=\sqc_{\D\in\co}\Si^\D$.
\nl
Note that $G$ acts naturally on $Z_\co$ (this action induces the conjugation action of $G$ on $\co$) and for any
$\D\in\co$ the obvious map $G^\D_{\ge0}\bsl\Si^\D\m G\bsl Z_\co$ is a bijection denoted by $\o_\D\lra\o$. (We 
use the fact the stabilizer in $G$
of an element $\D\in\co$ is equal to $G^\D_{\ge0}$.) Since $G^\D_{\ge0}\bsl\Si^\D$ is a finite set whose cardinal is a 
power of $2$ (see 1.9) we see that $G\bsl Z_\co$ is a finite set whose cardinal is a power of $2$. 
For any $\D\in D_G$ and $\o\in G\bsl Z_\co$ let $\fS^\D_{\o}$ be the inverse image of $\o_\D$ under the map
$\fS^\D@>>>\Si^\D$ in 2.3; we have a partition $\fS^\D=\sqc_{\o\in G\bsl Z_\co}\fS^\D_{\o}$. We set

$\cu_G^\co=\Ps_G(\sqc_{\D\in\co}\fS^\D)$,

$\cu_G^{\co,\o}=\Ps_G(\sqc_{\D\in\co}\fS^\D_{\o})$, ($\o\in G\bsl Z_\co$).
\nl 
The subsets $\cu_G^\co$ are called the {\it pieces } of $\cu_G$.
They form a partition of $\cu_G$ into subsets (which are unions of $G$-orbits) indexed by $\fA_G=\fA_{G'}$.
The subsets $\cu_{G.\o}^\co$ are called the {\it subpieces } of $\cu_G$.
When $\o$ varies in $G\bsl Z_\co$ the subpieces $\cu_{G,\o}^\co$ form a partition of a piece $\cu_G^\co$ into
subsets (which are unions of $G$-orbits); the number of these subsets is a power of $2$.

If $G$ is as in 1.3 then each piece of $\cu_G$ is a single $G$-orbit.
If $G$ is as in 1.4 ($p=2$) then each subpiece of $\cu_G$ is a single $G$-orbit..
If $G,V$ are as in 1.5 ($p=2$, $\dim V=9$) then there exist a subpiece of $\cu_G$ which is a union of two
$G$-orbits. 

In the remainder of this subsection we assume that 2.1(c) holds. 
Then for each $\co\in\fA_G$ and $n\ge1$, the piece $\cu_G^\co$ is $F$-stable and according to \cite{\LUU}, 
\cite{\LU}, the number $|(\cu_G^\co)^{F^n}|$ is a polynomial in $p^n$ with integer coefficients independent of 
$p,n$. From the definitions we have 

(a) $|(\cu_G^\co)^{F^n}|=|\co^{F^n}|\cdot|(G_{\ge3}^\D)^{F^n}|\cdot|(\Si^\D)^{F^n}|$ 
\nl
where $\D$ is any point of $\co^F$ (Note that $\co$ is $F$-stable by 2.1.)

\subhead 2.5\endsubhead  
Let $V,G,\fg$ be as in 1.3. A filtration of $V$ is a collection of subspaces $V_*=(V_{\ge a})_{a\in\ZZ}$ of $V$ 
such that $V_{\ge a+1}\sub V_{\ge a}$ for all $a$, $V_{\ge a}=0$ for some $a$, $V_{\ge a}=V$ for some $a$). If a 
filtration as above is given, we set $\gr_a(V_*)=V_{\ge a}/V_{\ge a+1}$ for any $a$ and 
$\gr(V_*)=\op_a\gr_a(V_*)$. Let $\fF(V)$ be the set of filtrations $V_*=(V_{\ge a})_{a\in\ZZ}$ of $V$ such that 
the grading $(\gr_a(V_*))$ of $\gr(V_*)$ is good (see 1.3) or equivalently such that there exists a good grading 
$(V_i)$ of $V$ with $V_{\ge a}=\op_{i;i\ge a}V_i$ for all $i$.

If $\d,\d'\in\fD_G$ correspond to the good gradings $(V_i)$, $(V'_i)$ of $V$ then $\d\si\d'$ if and only if
$\op_{i;i\ge a}V_i=\op_{i;i\ge a}V'_i$ for all $a\in\ZZ$. Setting $V_{\ge a}=\op_{i;i\ge a}V_i$ we see that
$(V_{\ge a})\in\fF(V)$ and that $(\si-\text{equivalence class of }\d)\m(V_{\ge a})$ is a bijection 
$D_G@>\si>>\fF(V)$.

For any $V_*=(V_{\ge a})\in\fF(V)$ let $\x(V_*)$ be the set of all $A\in\End(V)$ such that 
$A(V_{\ge a})\sub V_{\ge a+2}$ for any $a\in\ZZ$ and such that the map $\bA\in\End(\gr(V_*))_2$ induced by $A$ 
belongs to $\End(\gr(V_*))_2^0$. Note that $\x(V_*)\sub\cn_\fg$. Using 1.3(a) we see that in our case the 
following statement is equivalent to 2.4:

(a) {\it the map $\sqc_{V_*\in\fF(V)}\x(V_*)@>>>\cn_\fg$, $A\m A$, is a bijection.}
\nl
Let $A\m V_*^A$ be the map $\cn_\fg@>>>\fF(V)$ defined in \cite{\LUU, 2.3}. Then $A\m(A,V_*^A)$ is a well defined
map $\cn_\fg@>>>\sqc_{V_*\in\fF(V)}\x(V_*)$ which, by \cite{\LUU, 2.4} is an inverse of the map (a).

\subhead 2.6\endsubhead  
Let $V,(,),G,\fg=\fs(V)$ be as in 1.4. A filtration $V_*=(V_{\ge a})_{a\in\ZZ}$ of $V$ is said to be self-dual if 
$\{x\in V;(x,V_{\ge a})=0\}=V_{\ge1-a}$ for any $a$. If this is so, the associated vector space $\gr(V_*)$ has a 
unique symplectic form $(,)_0$ such that $(\gr_a(V_*),\gr_{a'}(V_*))_0=0$ if $a+a'\ne0$ and, for any 
$x\in\gr_a(V)$, $y\in\gr_{-a}(V)$, we have $(x,y)_0=(\dx,\dy)$ where $\dx\in V_{\ge a},\dy\in V_{\ge-a}$ are 
representatives of $x,y$. Note that $(,)_0$ is nondegenerate.
Let $\fF_s(V)$ be the set of self-dual filtrations $V_*=(V_{\ge a})_{a\in\ZZ}$ of $V$ such that 
the grading $(\gr_a(V_*))$ of $\gr(V_*)$ is $s$-good (see 1.4) or equivalently such that there exists an $s$-good
grading $(V_i)$ of $V$ with $V_{\ge a}=\op_{i;i\ge a}V_i$ for all $i$.

If $\d,\d'\in\fD_G$ correspond to the $s$-good gradings $(V_i)$, $(V'_i)$ of $V$ then $\d\si\d'$ if and only if
$\op_{i;i\ge a}V_i=\op_{i;i\ge a}V'_i$ for all $a\in\ZZ$. Setting $V_{\ge a}=\op_{i;i\ge a}V_i$ we see that
$(V_{\ge a})\in\fF_s(V)$ and that $(\si-\text{equivalence class of }\d)\m(V_{\ge a})$ is a bijection 
$D_G@>\si>>\fF_s(V)$.

Let $\cm(V)$ be the set of all nilpotent elements $A\in\End(V)$ such that $(Ax,y)+(x,Ay)+(Ax,Ay)=0$ for all $x,y$
in $V$ (or equivalently such that $1+A\in Sp(V)$). For any $V_*=(V_{\ge a})\in\fF_s(V)$ let $\ti\x(V_*)$ 
be the set of all $A\in\cm(V)$ such that $A(V_{\ge a})\sub V_{\ge a+2}$
for any $a\in\ZZ$ and such that the map $\bA\in\End(\gr(V_*))_2$ induced by $A$ belongs to $\End(\gr(V_*))_2^0$. 
Using 1.4(a) we see that in our case the following statement is equivalent to 2.4:

(a) {\it the map $\sqc_{V_*\in\fF_s(V)}\ti\x(V_*)@>>>\cm(V)$, $A\m A$, is a bijection.}
\nl
If $A\in\cm(V)$ then $V_*^A\in\fF(V)$ (see 2.5) is self-dual (see \cite{\LUU, 3.2(c)}) and the map
$\bA\in\End(\gr(V_*))_2$ induced by $A$ is in $\End(\gr(V_*))_2^0$ and is skew-adjoint with respect to $(,)_0$
(see \cite{\LUU, 3.2(d)}); this implies that $\dim\gr_a(V_*)$ is even when $a$ is even, hence $V_*^A\in\fF_s(V)$.
Then $A\m(A,V_*^A)$ is a well 
defined map $\cm(V)@>>>\sqc_{V_*\in\fF_s(V)}\ti\x(V_*)$ which, by \cite{\LUU, 2.4} is an inverse of the map (a).

\subhead 2.7\endsubhead  
Let $V,Q,(,),G=SO(V),\fg=\fo(V)$ be as in 1.5. A filtration $V_*=(V_{\ge a})_{a\in\ZZ}$ of $V$ is said to be a
$Q$-filtration if for any $a\ge1$ we have $Q|_{V_{\ge a}}=0$ and $\{x\in V;(x,V_{\ge a})=0\}=V_{\ge1-a}$. Then $Q$
induces (as in \cite{\LU, 1.5}) a nondegenerate quadratic form $\bQ:\gr V_*@>>>\kk$.

Let $\fF_o(V)$ be the set of all $Q$-filtrations $V_*=(V_{\ge a})_{a\in\ZZ}$ of $V$ such that 
the grading $(\gr_a(V_*))$ of $\gr(V_*)$ is $o$-good (see 1.5) or equivalently such that there exists an $o$-good
grading $(V_i)$ of $V$ with $V_{\ge a}=\op_{i;i\ge a}V_i$ for all $i$.

If $\d,\d'\in\fD_G$ correspond to the $o$-good gradings $(V_i)$, $(V'_i)$ of $V$ then $\d\si\d'$ if and only if
$\op_{i;i\ge a}V_i=\op_{i;i\ge a}V'_i$ for all $a\in\ZZ$. Setting $V_{\ge a}=\op_{i;i\ge a}V_i$ we see that
$(V_{\ge a})\in\fF_o(V)$ and that $(\si-\text{equivalence class of }\d)\m(V_{\ge a})$ is a bijection 
$D_G@>\si>>\fF_o(V)$. Let 
$$\cm(V)=\{A\in\End(V);A\text{ nilpotent}, 1+A\in SO(V)\}.$$
For any $V_*=(V_{\ge a})\in\fF_o(V)$ let $\et(V_*)$ be the set of all $A\in\cm(V)$ such that 
$A(V_{\ge a})\sub V_{\ge a+2}$ for any $a\in\ZZ$ and such that the map 
$\bA\in\End(\gr(V_*))_2$ induced by $A$ belongs to $\fo(\gr(V_*))_2^0$ (the last set is defined in terms of the
grading $(\gr_a(V_*))$ and the quadratic form on $\gr(V_*)$ induced by $Q$). Using 1.3(a) we see that in our case
the following statement is equivalent to 2.4:

(a) {\it the map $\sqc_{V_*\in\fF_o(V)}\et(V_*)@>>>\cm(V)$, $A\m A$, is a bijection.}
\nl
Assume first that $p\ne2$. If $A\in\cm(V)$ then $V_*^A\in\fF(V)$ (see 2.5) is in $\fF_o(V)$ (see 
\cite{\LU, 3.3}). Then $A\m(A,V_*^A)$ is a well defined map $\cm(V)@>>>\sqc_{V_*\in\fF_o(V)}\et(V_*)$ which, 
by \cite{\LUU, 2.4} is an inverse of the map (a). 

Next we assume that $p=2$. If $A\in\cm(V)$ we define a filtration $V^*_A=(V^{\ge a}_A)$ of $V$ as in
\cite{\LU, 2.5}. (Note that in general $V^*_A$ is not equal to $V_*^A$ of 2.5.) From \cite{\LU, 2.5(a)} 
we see that $V^*_A$ is a $Q$-filtration and from \cite{\LU, 2.9} we see that $V_*^A\in\fF_o(V)$ and
$A\in\et(V_*^A)$. Then 
$A\m(A,V^*_A)$ is a well defined map $\cm(V)@>>>\sqc_{V_*\in\fF_o(V)}\et(V_*)$ which, by \cite{\LU, 2.10} is an
inverse of the map (a). (An alternative proof of (a) is given in 2.9-2.10.)

\subhead 2.8\endsubhead  
Let $V,Q,(,),G=SO(V)$ be as in 1.5. In the remainder of this section we assume that $p=2$.

For any nonzero nilpotent element $T\in\End(V)$ let $e=e_T$ be the smallest integer 
$\ge2$ such that $T^e=0$; let $f=f_T$ be the smallest integer $\ge1$ s.t. $QT^f=0$. 
We associate to $T$ a subset $H_T$ of $V$ as follows.

(i) $H_T=\{x\in V;T^{e-1}x=0\}$ if $e\ge 2f$,

(ii) $H_T=\{x\in V;T^{e-1}x=0,QT^{f-1}x=0\}$ if $e=2f-1$,

(iii) $H_T=\{x\in V;QT^{f-1}x=0\}$ if $e<2f-1$.
\nl
We fix an $o$-good grading $V=\op_iV_i$ such that $V\ne V_0$. For $a\in\ZZ$ let 
$V_{\ge a}=V_a+V_{a+1}+\do$. Let $m\ge1$ be the largest integer such that
$V_m\ne0$. Let $A\in\fo(V)_2^0$. Note that $A\ne0$ and $A$ is nilpotent hence $\be:=e_A,\baf:=f_A,H_A$ are 
defined.

From the definition of $\fo(V)_2^0$ we see that:

(a) {\it If $m$ is odd then $\be=2\baf=m+1$. If $m$ is even then either $\be=2\baf-1=m+1$ or $\be<2\baf-1=m+1$. In 
any case, $V_{\ge-m+1}=H_A$.}
\nl
Let $T\in\End(V)$ be such that $C:=T-A$ satisfies $C(V_i)\sub V_{\ge i+3}$ for any $i$. Then $T\ne0$, 
$T(V_i)\sub V_{\ge i+2}$ for any $i$ hence $T$ is nilpotent so that $e=e_T,f=f_T,H_T$ are defined. We shall prove
successively the statements (b)-(j).

(b) {\it If $2n>m$ then $QT^n=0$, $QA^n=0$. Hence if $m$ is odd then $f\le(m+1)/2$, $\baf\le(m+1)/2$; if $m$ is 
even then $f\le(m+2)/2$, $bof\le(m+2)/2$.}
\nl
To show that $QT^n=0$ it is enough to show that $Q(T^nx)=0$ whenever $x\in V_i$, $i\ge-m$ and $(T^nx,T^nx')=0$ 
whenever $x\in V_i,x'\in V_j$, $i,j\ge-m$. This follows from $T^nx\in V_{\ge1}$ and $Q|_{V_{\ge1}}=0$, 
$(V_{\ge1},V_{\ge1})=0$. The equality $QA^n=0$ is proved in the same way.

The following result is immediate.

(c) {\it Let $P_n\in\End(V)$ be a product of $n$ factors of which at least one is $C$ and the remaining ones are 
$A$. Then $P_n(V_i)\sub V_{\ge i+2n+1}$ for all $i$. Hence if $n\ge m$ then $P_n=0$.}

(d) {\it If $n\ge m$ then $T^n=A^n$.}
\nl
Indeed, $T^n=A^n+\sum P_n$ where each term $P_n$ is as in (c). Hence the result follows from (c).

(e) {\it If $\be\ge2\baf-1$ then $T^{\be-1}=A^{\be-1}\ne0$, $T^{\be}=A^{\be}=0$; hence $e=\be$. Moreover if $\be=2\baf$ 
then $H_T=H_A$.}
\nl 
In our case we have $\be=m+1$ (see (a)) hence the first sentence follows from (d). If $\be=2\baf$ then $m$ is odd
and $e=\be=2\baf=m+1\ge 2f$ (see (b)) hence $e\ge2f$ and 
$$\align&H_T=\{x\in V;T^{e-1}x=0\}=\{x\in V;T^{\be-1}x=0\}\\&=\{x\in V;A^{\be-1}x=0\}=H_A.\endalign$$
(f) {\it If $2n\ge m$, $P_n$ is as in (c) and $P'_n$ is like $P_n$ then $QP_n(V)=0$, $(P_n(V),P'_n(V))=0$, 
$(A^n(V),P_n(V))=0$.}
\nl
By (c) we have $P_nV_i\sub V_{\ge1}$ if $i\ge-m$ hence $P_n(V)\sub V_{\ge1}$. Also, $A^nV\sub V_{\ge0}$. It is 
then enough to use: $Q|_{V_{\ge1}}=0$, $(V_{\ge1},V_{\ge0})=0$.

(g) {\it If $2n\ge m$ then $QT^n=QA^n$.}
\nl
We have $T^n=A^n+\sum P_n$ where each term $P_n$ is as in (c). Hence for $x\in V$ we have
$QT^nx=QA^nx+\sum QP_nx+\sum(A^n,P_nx)+\sum(P_nx,P'_nx)$ with $P'_n$ like $P_n$ and we use (f).

(h) {\it If $\be<2\baf$ then $QT^{\baf-1}=QA^{\baf-1}\ne0$, $QT^{\baf}=QA^{\baf}=0$; hence $f=\baf$.}
\nl
In this case we have $2(\baf-1)=m$, $2\baf>m$. Hence result follows from (g).

(i) {\it If $\be<2\baf-1$ then $e<2f-1$.}
\nl
In this case we have $f=\baf=(m+2)/2$. We must show that $e<m+1$. 
By (d) we have $T^m=A^m$. Since $\be\le m$ (see (a)) we have $A^m=0$ hence $T^m=0$ and $e<m+1$.

(j) {\it If $2\baf=\be$ then $e=\be=m+1\ge f/2$. If $2\baf-1=\be$ then $e=\be=m+1$, $f=\baf$ hence $2f-1=e$. 
If $2\baf-1>\be$ then $f=\baf=(m+2)/2$ and $2f-1>e$. In each case we have $H_T=H_A$ hence $H_T=V_{\ge-m+1}$ 
and $m=\max(e-1,2f-2)$.}

\subhead 2.9\endsubhead
We give an alternative proof of the injectivity of the map 2.7(a). 
We argue by induction on $\dim V$. If $\dim V\le1$, the result is 
trivial. Assume now that $\dim V\ge2$. Let $T\in\cm(V)$ and let $V_*=(V_{\ge a})$, $\tV_*=(\tV_{\ge a})$ be two 
filtrations in $\fF_o(V)$ such that $T\in\et(V_*)$ and $T\in\et(\tV_*)$. We must show that $V_*=\tV_*$. Let 
$\bT\in\End(\gr V_*)_2$, $\bT_1\in\End(\gr\tV_*)_2$ be the endomorphisms induced by $T$. If $\bT=0$ then, using 
the fact that $T\in\et'(V_*)$ we see that $\gr_a(V_*)=0$ for $a\ne0$ hence $V_{\ge1}=0,V_{\ge 0}=V$; since 
$TV=T(V_{\ge 0})\sub V_{\ge 2}=0$ we see that $T=0$, hence $\bT_1=0$ and $\tV_{\ge1}=0,\tV_{\ge 0}=V$; thus
$V_*=\tV_*$ as desired. Similarly, if $\bT_1=0$ then $V_*=\tV_*$ as desired. Thus we can assume
that $\bT\ne0$, $\bT_1\ne0$. Hence $\gr_a(V_*)\ne0$ for some $a\ne0$ and $\gr_a(\tV_*)\ne0$ for some $a\ne0$. Let
$m\ge1$ be the largest integer such that $\gr_m(V_*)\ne0$. Let $\tm\ge1$ be the largest integer such that 
$\gr_{\tm}(\tV_*)\ne0$. Using 2.8(j) we see that $V_{\ge-m+1}=H_T=\tV_{\ge-\tm+1}$, $m=\max(e_T-1,2f_T-2)=\tm$.
It follows that $m=\tm$ and $V_{\ge-m+1}=\tV_{\ge-m+1}$. Since $V_*$ is a $Q$-filtration
we have $V_{\ge m}=\{x\in V;(x,V_{\ge-m+1})=0;Q(x)=0\}$, see \cite{\LU, 1.4(b)}. Similarly we have
$\tV_{\ge m}=\{x\in V;(x,\tV_{\ge-m+1})=0;Q(x)=0\}$. Hence $V_{\ge m}=\tV_{\ge m}$. Let 
$V'=V_{\ge-m+1}/V_{\ge m}=\tV_{\ge-m+1}/\tV_{\ge m}$. Note that $V'$ has a natural nondegenerate quadratic form 
induced by $Q$. We set $V'_{\ge a}=\text{image of $V_{\ge a}$ under }V_{\ge-m+1}@>>>V'$ (if $a\ge-m+1$), 
$V'_{\ge a}=0$ (if $a<-m+1$). We set $\tV'_{\ge a}=\text{image of $\tV_{\ge a}$ under }\tV_{\ge-m+1}@>>>V'$ (if 
$a\ge-m+1$), $\tV'_{\ge a}=0$ (if $a<-m+1$). Then $V'_*=(V'_{\ge a})$, $\tV'_*=(\tV'_{\ge a})$ are filtrations in
$\fF_o(V')$. Also $T$ induces an element $T'\in\cm(V')$ and we have $T'\in\et(V'_*)$, $T'\in\et(\tV'_*)$. Note 
also that $\dim V'<\dim V$. By the induction hypothesis we have $V'_*=\tV'_*$. It follows that 
$V_{\ge a}=\tV_{\ge a}$ for any $a\ge-m+1$. If $a<-m+1$ we have $V_{\ge a}=\tV_{\ge a}=V$. Hence $V_*=\tV_*$, as
desired. Thus the map 2.7(a) is injective. 

\subhead 2.10\endsubhead
We give an alternative proof of the surjectivity of the map 2.7(a). 
By a standard argument we can assume that $\kk$ is an algebraic
closure of the field $\FF_2$ with $2$ elements. We can also assume that $\dim V\ge2$. We choose an 
$\FF_2$-rational structure on $V$ such that $Q$ is defined and split over $\FF_2$. Then the Frobenius map
relative to the $\FF_2$-structure acts naturally and compatibly on the source and target of the map 2.7(a). 
We denote each of these actions by $F$. It is enough to show that for any $n\ge1$ the map
$\a_n:(\sqc_{V_*\in\fF_o(V)}\et(V_*))^{F^n}@>>>\cm(V)^{F^n}$, $A\m A$ is a bijection. (Here $()^{F^n}$ is the 
fixed point set of $F^n$.) Since $\a_n$ is injective (see 2.9) it is enough to show that
$|(\sqc_{V_*\in\fF_o(V)}\et(V_*))^{F^n}|=|\cm(V)^{F^n}|$.
By \cite{\ST, 15.1} we have $|\cm(V)^{F^n}|=2^{nr}$ where $r$ is the number of roots of $G$. It is enough to show
that 

(a) $\sum_{V_*\in\fF_o(V)^{F^n}}|\et(V_*)^{F^n}|=2^{nr}$. 
\nl
Now the left hand side of (a) makes sense when $\kk$ is replaced by an algebraic closure of the prime field with 
$p'$ elements where $p'$ is any prime number; moreover this more general expression is a polynomial in $p'{}^n$
with integer coefficients independent of $p',n$ (see \cite{\LU, 1.5(e)}). When $p'\ne2$, this more general
expression is equal to $p'{}^{nr}$ since 2.7(a) is already known to be a bijection in this case. But then (a) 
follows from this more general equality by specializing $p'{}^n$ (viewed as indeterminate) to $2^n$.

\subhead 2.11\endsubhead
Assume that $G$ is as in 2.4 and that $\D,\ti\D\in D_G$ are conjugate under $G$. We show:

(a) {\it if $u\in\fS^\D$, $u\in G_{\ge2}^{\ti\D}$ for some $u$ then $\D=\ti\D$.}
\nl
We can assume that $G$ is as in 1.3, 1.4 or 1.5. 

Assume first that $G,\fg,V$ are as in 1.3. We prove (a) by induction on $\dim V$. If $V=0$ the result is trivial.
We 
now assume that $V\ne0$. Let $V_*=(V_{\ge a})$, $\tV_*=(\tV_{\ge a})$ be the objects of $\fF(V)$ corresponding to 
$\D,\ti\D$. Let $m$ be the largest integer $\ge0$ such that $\gr_mV_*\ne0$. If $m=0$ then $V_a=V$ for $a\le m$ and
$V_a=0$ for $a>m$. Since $\tV_*$, $V_*$ are $G$-conjugate we have $\tV_*=V_*$. Thus we may assume that $m\ge1$. We
have $V_{\ge-m}=V$, $V_{\ge m+1}=0$. Since $\tV_*$, $V_*$ are $G$-conjugate we see that $\tV_{\ge-m}=V$, 
$\tV_{\ge m+1}=0$. We set $x=u-1$. From $x\in\x(V_*)$ we deduce $V_{\ge-m+1}=\ker(x^m:V@>>>V)$, $V_{\ge m}=x^mV$.
From $x(\tV_{\ge a})\sub\tV_{\ge a+2}$ for all $a$ we deduce 
$x^m(\tV_{\ge-m+1})\sub\tV_{\ge-m+1+2m}=0$, $x^m(V)=x^m(\tV_{\ge-m})\sub\tV_{\ge-m+2m}=\tV_{\ge m}$. Thus, 
$\tV_{\ge-m+1}\sub\ker(x^m:V@>>>V)$ and $\tV_{\ge-m+1}\sub V_{\ge-m+1}$; moreover, $V_{\ge m}\sub\tV_{\ge m}$. 
Since $\tV_*$, $V_*$ are $G$-conjugate we have $\dim\tV_{\ge-m+1}=\dim V_{\ge-m+1}$, 
$\dim\tV_{\ge m}=\dim V_{\ge m}$ hence $\tV_{\ge-m+1}=V_{\ge-m+1}$, $\tV_{\ge m}=V_{\ge m}$. Let 
$V'=V_{\ge-m+1}/V_{\ge m}=\tV_{\ge-m+1}/\tV_{\ge m}$. Now $V_*,\tV_*$ give rise in an obvious way to two elements
$V'_*,\tV'_*$ of $\fF(V')$ and $u$ gives rise to a unipotent element $u'\in GL(V')$ such that $u'-1\in\x(V'_*)$ 
and
$(u'-1)(\tV'_{\ge a})\sub\tV'_{\ge a+2}$ for all $a$. Since $\dim V'<\dim V$, the induction hypothesis shows that 
$V'_*=\tV'_*$. It follows that $V_*=\tV_*$. This proves (a) in our case. 

Assume next that $G,\fg,V,(,)$ are as in 1.4. In this case the proof is essentially the same as in the case
of 1.3. The same applies in the case where $G,\fg,V,Q,(,)$ are as in 1.5 and $p\ne2$.

In the remainder of this subsection we assume that $G,V,Q,(,)$ are as in 1.5 and $p=2$.
We prove (a) by induction on $\dim V$. If $V=0$ the result is trivial. We now assume that $V\ne0$. Let 
$V_*=(V_{\ge a})$, $\tV_*=(\tV_{\ge a})$ be the objects of $\fF_o(V)$ corresponding to 
$\D,\ti\D$. Let $m$ be the largest integer $\ge0$ such that $\gr_mV_*\ne0$. If $m=0$ then $V_a=V$ for $a\le m$ and
$V_a=0$ for $a>m$. Since $\tV_*$, $V_*$ are $G$-conjugate we have $\tV_*=V_*$. Thus we may assume that $m\ge1$. 
We have $V_{\ge-m}=V$. Since $\tV_*$, $V_*$ are $G$-conjugate we see that $\tV_{\ge-m}=V$. Let $x=u-1$. From 
$x(\tV_{\ge a})\sub\tV_{\ge a+2}$ for all $a$ we deduce $x^m(\tV_{\ge-m+1})\sub\tV_{\ge-m+1+2m}=0$; moreover if 
$m$ is even we have 
$Qx^{m/2}(\tV_{\ge-m+1})\sub Q\tV_{\ge-m+1+m}=Q\tV_{\ge1}=0$.
\nl
Thus,
$\tV_{\ge-m+1}\sub\{v\in V;x^mv=0\}$ and if $m$ is even, $\tV_{\ge-m+1}\sub\{v\in V;Qx^{m/2}v=0\}$. Recall from 
2.8(j) that $m=\max(e_x-1,2f_x-2)$. If $e_x\ge 2f_x$ we have $m=e_x-1$ and

$\tV_{\ge-m+1}\sub\{v\in V;x^{e_x-1}v=0\}=H_x=V_{\ge-m+1}$.
\nl
If $e_x<2f_x-1$ we have $m=2f_x-2$ and

$\tV_{\ge-m+1}\sub\{v\in V;Qx^{f_x-1}v=0\}=H_x=V_{\ge-m+1}$.
\nl
If $e_x=2f_x-1$ we have $m=2f_x-2=e_x-1$ and

$\tV_{\ge-m+1}\sub\{v\in V;x^{e_x-1}v=0\}\cap\{v\in V;Qx^{f_x-1}v=0\}=H_x=V_{\ge-m+1}$.
\nl
Thus in any case we have $\tV_{\ge-m+1}\sub V_{\ge-m+1}$. It follows that $\tV_{\ge m}=V_{\ge m}$.
Let $V'=V_{\ge-m+1}/V_{\ge m}=\tV_{\ge-m+1}/\tV_{\ge m}$. Note that $Q$ induces naturally a nondegenerate
quadratic form on $V'$. Also $V_*,\tV_*$ give rise in an obvious way to two elements $V'_*,\tV'_*$ of 
$\fF_o(V')$ and $u$ gives rise to a unipotent element $u'\in SO(V')$ such that $(u'-1)\in\et'(V'_*)$ and 
$(u'-1)(\tV'_{\ge a})\sub\tV'_{\ge a+2}$ for all $a$. Since $\dim V'<\dim V$, the induction hypothesis shows that 
$V'_*=\tV'_*$. It follows that $V_*=\tV_*$. This proves (a) in our case.

\head Appendix: The pieces in the nilpotent variety of $\fg$\endhead
{\it by G. Lusztig and T. Xue}

\subhead A.1\endsubhead
For any $\D\in D_G$ let $\fs^\D\sub\fg_{\ge2}^\D$ be the inverse image of $\Si^\D$ (see 2.3) under the obvious map 
$\fg_{\ge2}^\D@>>>\fg_{\ge2}^\D/\fg_{\ge3}^\D$. Now $\fs^\D$ is stable under the $\Ad$-action of $G^\D_{\ge0}$ on 
$\fg_{\ge2}^\D$ and $x\m x$ is a map
$$\Ps_\fg:\sqc_{\D\in D_G}\fs^\D@>>>\cn_\fg.$$

\proclaim{Theorem A.2} Assume that $\tG^{der}$ (see 1.1) is a product of almost simple groups of type $A,B,C,D$.
Then $\Ps_\fg$ is a bijection.
\endproclaim
The general case reduces easily to the case where $G$ is almost simple of type $A,B,C$ or $D$. Moreover
we can assume that $G$ is one of the groups $GL(V)$, $Sp(V)$, $SO(V)$ 
in 1.3-1.5. The proof in these cases will be given
in A.3-A.4. We expect that the theorem holds without restriction on $G$.

\subhead A.3\endsubhead
If $V,G,\fg$ are as in 1.3, the proof of A.2 is exactly as in 2.5. Now assume that $V,(,),G,\fg=\fs(V)$ are as in 1.4. 
For any $V_*=(V_{\ge a})\in\fF_s(V)$ let $\ti\x'(V_*)$ be the set of all $A\in\fs(V)$ such that 
$A(V_{\ge a})\sub V_{\ge a+2}$
for any $a\in\ZZ$ and such that the map $\bA\in\End(\gr(V_*))_2$ induced by $A$ belongs to $\End(\gr(V_*))_2^0$. 
Using 1.4(a) we see that in our case the following statement is equivalent to A.2:

(a) {\it the map $\sqc_{V_*\in\fF_s(V)}\ti\x'(V_*)@>>>\cn_\fg$, $A\m A$, is a bijection;}
\nl
If $A\in\cn_\fg$ then $V_*^A\in\fF(V)$ (see 2.5) is self-dual and the map $\bA\in\End(\gr(V_*))_2$ induced by $A$
 is in $\End(\gr(V_*))_2^0$ and is skew-adjoint with respect to $(,)_0$ (the proofs are completely similar to 
those in \cite{\LUU, 3.2(c)}, \cite{\LUU, 3.2(d)}; this implies that $\dim\gr_a(V_*)$ is even when $a$ is even, 
hence $V_*^A\in\fF_s(V)$. Then $A\m(A,V_*^A)$ is a well 
defined map $\cn_\fg@>>>\sqc_{V_*\in\fF_s(V)}\ti\x'(V_*)$ which, by \cite{\LUU, 2.4} is an inverse of the map (a).

\subhead A.4\endsubhead  
Let $V,Q,(,),G=SO(V),\fg=\fo(V)$ be as in 1.5. 
For any $V_*=(V_{\ge a})\in\fF_o(V)$ let $\et'(V_*)$ be the set of all $A\in\fo(V)$ such that 
$A(V_{\ge a})\sub V_{\ge a+2}$ for any $a\in\ZZ$ and such that the map 
$\bA\in\End(\gr(V_*))_2$ induced by $A$ belongs to $\fo(\gr(V_*))_2^0$ (the last set is defined in terms of the
grading $(\gr_a(V_*))$ and the quadratic form on $\gr(V_*)$ induced by $Q$). Using 1.5(a) we see that in our case
the following statements are equivalent to A.2:

(a) {\it the map $\sqc_{V_*\in\fF_o(V)}\et'(V_*)@>>>\cn_\fg$, $A\m A$, is a bijection.}
\nl
Assume first that $p\ne2$. If $A\in\cn_\fg$ then $V_*^A\in\fF(V)$ (see 2.5) is in $\fF_o(V)$ (by an argument
entirely similar to that in \cite{\LU, 3.3}). Then $A\m(A,V_*^A)$ is a well defined map 
$\cn_\fg@>>>\sqc_{V_*\in\fF_o(V)}\et'(V_*)$ which, by \cite{\LUU, 2.4} is an inverse of the map (a).
It remains to prove (a) in the case where $p=2$. 

We show that the map (a) is injective. The proof is almost the same as that in 2.9. We argue by induction on 
$\dim V$. If $\dim V\le1$ the result is trivial. Assume now that $\dim V\ge2$. Let $T\in\fo(V)$ and let 
$V_*=(V_{\ge a})$, $\tV_*=(\tV_{\ge a})$ be two filtrations in $\fF_o(V)$ such that $T\in\et'(V_*)$ and 
$T\in\et'(\tV_*)$. We must show that $V_*=\tV_*$. Let $\bT\in\End(\gr V_*)_2$, $\bT_1\in\End(\gr\tV_*)_2$ be the 
endomorphisms induced by $T$. If $\bT=0$ or $\bT'=0$ then, as in the proof in 2.9 we see that $V_*=\tV_*$, as 
desired. Thus we can assume that $\bT\ne0$, $\bT_1\ne0$. Hence $\gr_a(V_*)\ne0$ for some $a\ne0$ and 
$\gr_a(\tV_*)\ne0$ for some $a\ne0$. Let $m\ge1$ be the largest integer such that $\gr_m(V_*)\ne0$. Let $\tm\ge1$
be the largest integer such that $\gr_{\tm}(\tV_*)\ne0$. Using 2.8(j) we see that 
$V_{\ge-m+1}=H_T=\tV_{\ge-\tm+1}$, $m=\max(e_T-1,2f_T-2)=\tm$.
It follows that $m=\tm$ and $V_{\ge-m+1}=\tV_{\ge-m+1}$. Since $V_*$ is a $Q$-filtration
we have $V_{\ge m}=\{x\in V;(x,V_{\ge-m+1})=0;Q(x)=0\}$, see \cite{\LU, 1.4(b)}. Similarly we have
$\tV_{\ge m}=\{x\in V;(x,\tV_{\ge-m+1})=0;Q(x)=0\}$. Hence $V_{\ge m}=\tV_{\ge m}$. Let 
$V'=V_{\ge-m+1}/V_{\ge m}=\tV_{\ge-m+1}/\tV_{\ge m}$. Note that $V'$ has a natural nondegenerate quadratic form 
induced by $Q$. We set $V'_{\ge a}=\text{image of $V_{\ge a}$ under }V_{\ge-m+1}@>>>V'$ (if $a\ge-m+1$), 
$V'_{\ge a}=0$ (if $a<-m+1$). We set $\tV'_{\ge a}=\text{image of $\tV_{\ge a}$ under }\tV_{\ge-m+1}@>>>V'$ (if 
$a\ge-m+1$), $\tV'_{\ge a}=0$ (if $a<-m+1$). Then $V'_*=(V'_{\ge a})$, $\tV'_*=(\tV'_{\ge a})$ are filtrations in
$\fF_o(V')$. Also $T$ induces an element $T'\in\fo(V')$ and we have $T'\in\et'(V'_*)$, $T'\in\et'(\tV'_*)$. Note 
also that $\dim V'<\dim V$. By the induction hypothesis we have $V'_*=\tV'_*$. It follows that 
$V_{\ge a}=\tV_{\ge a}$ for any $a\ge-m+1$. If $a<-m+1$ we have $V_{\ge a}=\tV_{\ge a}=V$. Hence $V_*=\tV_*$, as
desired. Thus the map (a) is injective. 

We show that the map (a) is surjective. The proof is somewhat similar to that in 2.10. By a standard argument we 
can assume that $\kk$ is an algebraic closure of the field $\FF_2$ with $2$ elements. We can also assume that 
$\dim V\ge2$. We choose an $\FF_2$-rational structure on $V$ such that $Q$ is defined and split over $\FF_2$. Then
the Frobenius map relative to the $\FF_2$-structure acts naturally and compatibly on the source and target of the
maps (a) and 2.7(a). We denote each of these actions by $F$. It is enough to show that for any $n\ge1$ the map
$(\sqc_{V_*\in\fF_o(V)}\et'(V_*))^{F^n}@>>>\cn_\fg^{F^n}$, $A\m A$ is a bijection. (Here $()^{F^n}$ is the fixed
point set of $F^n$.) Since the last map is injective (by the previous paragraph) it is enough to show that
$|\sqc_{V_*\in\fF_o(V)}\et'(V_*))^{F^n}|=|\cn_\fg^{F^n}|$. From 2.10 we see that 
$|\sqc_{V_*\in\fF_o(V)}\et(V_*)^{F^n}|=2^{nr}$ ($r$ as in 2.10). According to \cite{\SPR} we have also 
$|\cn_\fg^{F^n}|=2^{nr}$. It is then enough to show that

$|\sqc_{V_*\in\fF_o(V)}\et'(V_*))^{F^n}|=|\sqc_{V_*\in\fF_o(V)}\et(V_*)^{F^n}|$
\nl
or equivalently that

$\sum_{V_*\in\fF_o(V)^{F^n}}|\et'(V_*)^{F^n}|=\sum_{V_*\in\fF_o(V)^{F^n}}|\et(V_*)^{F^n}|$.
\nl
It is enough to show that for any $V_*\in\fF_o(V)^{F^n}$ we have $|\et'(V_*)^{F^n}|=|\et(V_*)^{F^n}|$. From the 
definitions we have $|\et(V_*)^{F^n}|=|\fo(\gr(V_*)_2^0|2^{nd}$ where $d$ is the dimension of the vector space 
$\{C\in\fo(V);C(V_{\ge a})\sub V_{\ge a+3}\qua\frl a\}$. From \cite{\LU, 1.5(c),(d)} we see that 
$|\et(V_*)^{F^n}|=|\fo(\gr(V_*)_2^0|2^{nd}$ where $d$ is as above. We see that 
$|\et'(V_*)^{F^n}|=|\et(V_*)^{F^n}|$. This proves the surjectivity of the map (a) and completes the proof of A.2.

\subhead A.5\endsubhead
Assume that $G$ is as in 2.4 and that $\D,\ti\D\in D_G$ are conjugate under $G$. The following statement is a 
strengthening of the statement that $\Ps_\fg$ in A.2 is injective:

(a) {\it  If $x\in\fs^\D$, $x\in\fg_{\ge2}^{\ti\D}$ for some $x$ then $\D=\ti\D$.}
\nl
The proof is essentially the same as that of 2.11(a).

\subhead A.6\endsubhead
Let $G$ be as in A.2. Let $\co\in\fA_G$ (see 2.3). For any $\D\in D_G$ and $\o\in G\bsl Z_\co$ (see 2.3) let 
$\fs^\D_{\o}$ be the inverse image of $\o_\D$ (see 2.3) under the map $\fs^\D@>>>\Si^\D$ in A.1; we have a 
partition $\fs^\D=\sqc_{\o\in G\bsl Z_\co}\fs^\D_{\o}$. We set

$\cn_\fg^\co=\Ps_\fg(\sqc_{\D\in\co}\fs^\D)$,

$\cn_\fg^{\co,\o}=\Ps_G(\sqc_{\D\in\co}\fs^\D_{\o})$, ($\o\in G\bsl Z_\co$).
\nl 
The subsets $\cn_\fg^\co$ are called the {\it pieces } of $\cn_\fg$. They form a partition of $\cn_\fg$ into 
subsets (which are unions of $G$-orbits) indexed by $\fA_G=\fA_{G'}$. The subsets $\cn_{\fg,\o}^\co$ are called 
the {\it subpieces } of $\cn_\fg$. They are in bijection with the subpieces of $G$. When $\o$ varies in 
$G\bsl Z_\co$ the subpieces $\cn_{\fg,\o}^\co$ form a partition of a piece $\cn_\fg^\co$ into subsets (which are 
unions of $G$-orbits); the number of these subsets is a power of $2$.

In the remainder of this subsection we assume that 2.1(c) holds. Then for each $\co\in\fA_G$ and $n\ge1$, the 
piece $\cn_\fg^\co$ is $F$-stable and from the definitions we have 
$|\cn_\fg^\co)^{F^n}|=|\co^{F^n}|\cdot|(\fg_{\ge3}^\D)^{F^n}|\cdot|(\Si^\D)^{F^n}|$ where $\D$ is any point of 
$\co^F$. (Note that $\co$ is $F$-stable by 2.1.) We have $|(\fg_{\ge3}^\D)^{F^n}|=|G_{\ge3}^\D)^{F^n}|$. Using 
this and 2.4(a) we see that $|(\cn_\fg^\co)^{F^n}|=|(\cu_G^\co)^{F^n}|$. In particular, 
$|(\cn_\fg^\co)^{F^n}|$ is a polynomial in $p^n$ with integer coefficients independent of $p,n$.

\enddocument